\documentclass{amsart}

\usepackage{amsmath}
\usepackage{paralist}
\usepackage{amsfonts}
\usepackage{MnSymbol}
\usepackage{graphics} 
\usepackage{epsfig} 
\usepackage{psfrag}
\usepackage{pdfsync}
\usepackage[usenames]{color}
\usepackage{cancel}
\usepackage{mathrsfs}
\DeclareMathAlphabet{\mathpzc}{OT1}{pzc}{m}{it}

\definecolor{arancio}{rgb}{0.90,0.50,0.20}

\definecolor{blu}{rgb}{0.,0.,1.}

\definecolor{pavone}{rgb}{0.00,0.00,0.63}

\definecolor{malva}{rgb}{0.10,0.50,0.50}

\definecolor{rosso}{rgb}{1.,0.,0.}
\newcommand{\rosso}[1]{{\color{red}{#1}}}

\definecolor{geranio}{rgb}{0.90,0.00,0.20}

\definecolor{cerulean}{rgb}
{0.0, 0.48, 0.65}

\newtheorem{theorem}{Theorem}[section]

\newtheorem{lem}[theorem]{Lemma}
\newtheorem{prop}[theorem]{Proposition}

\theoremstyle{definition}
\newtheorem{definition}{Definition}[section]
\newtheorem{remark}{Remark}[section]

\newcommand{\N}{\mathbb{N}}
\newcommand{\R}{\mathbb{R}}

\date{\today}

\newcommand{\bcl}{\begin{center}}
\newcommand{\ecl}{\end{center}}
\newcommand{\brl}{\begin{right}}
\newcommand{\erl}{\end{right}}
\newcommand{\ben}{\begin{enumerate}}
\newcommand{\barr}{\begin{array}}
\newcommand{\earr}{\end{array}}
\newcommand{\btab}{\begin{tabular}}
\newcommand{\etab}{\end{tabular}}
\newcommand{\bdoc}{\begin{document}}
\newcommand{\edoc}{\end{document}}
\newcommand{\beqy}{\begin{eqnarray}}

\newcommand{\beq}{\begin{equation}}
\newcommand{\beqi}{\begin{eqnarray*}}
\newcommand{\bitem}{\begin{itemize}}
\newcommand{\brem}{\begin{remark}}
\newcommand{\erem}{\end{remark}}
\newcommand{\eitem}{\end{itemize}}
\newcommand{\nln}{\newline}
\newcommand{\newt}{\newtheorem}

\renewcommand{\a }{\alpha }
\renewcommand{\b }{\beta }
\newcommand{\g }{\gamma}
\newcommand{\G }{\Gamma }
\renewcommand{\d }{\delta }
\newcommand{\D }{\Delta }
\newcommand{\e }{\epsilon }
\newcommand{\z }{\zeta }
\renewcommand{\l }{\lambda  }
\renewcommand{\L }{\Lambda }
\newcommand{\m }{\mu }
\newcommand{\n }{\nu }
\renewcommand{\r }{\rho }
\newcommand{\s }{\sigma }
\newcommand{\Sig }{\Sigma }
\renewcommand{\t }{\tau }
\newcommand{\vT}{\Theta }
\renewcommand{\o }{\omega }
\renewcommand{\O }{\Omega }
\newcommand{\pa }{\partial  }

\newcommand{\supp}{\text{\rm supp}\,}
\newcommand{\sgn}{\text{\rm sgn}\,}

\newcommand{\red}[1]{{\color{red}{#1}}}
\newcommand{\blue}[1]{{\color{blue}{#1}}}

\title[Monotonicity results] {Monotonicity results \\ for semilinear parabolic equations on metric graphs}

\author[Punzo]{Fabio Punzo}
\address{Dipartimento di Matematica, Politecnico di Milano \\ Via Bonardi 9, I-20133 Milano, Italy}
\email{fabio.punzo@polimi.it}
\author[Tesei]{Alberto Tesei}
\address{Dipartimento di Matematica ``G. Castelnuovo", Universit\`a Sapienza di Roma\\ P.le A. Moro 5, I-00185 Roma, Italy, and
Istituto per le Applicazioni del Calcolo ``M. Picone", CNR,
Via dei Taurini 19, I-00185 Roma, Italy}
\email{alberto.tesei@lincei.it}



\date{\today}

\begin{document}

\bibliographystyle{h-elsevier2}

\begin{abstract} We prove monotonicity results for semilinear parabolic problems on locally finite connected metric graphs. Applications to regular metric trees are discussed.
\end{abstract}

\maketitle


\section{introduction}

 In this note we prove monotonicity results for the semilinear parabolic problem
\begin{subequations}\label{CN}
\begin{equation}\label{diffeq}
u_t  = \Delta u+ f(u) \quad\textrm{in}\,\,\mathcal{G}\times \R_+
\,,
\end{equation}
\begin{equation}\label{inco}
u \, = u_0 \quad \textrm{in}\,\, \mathcal{G}\times \{0\}\,.
\end{equation}
\end{subequations}
Here $\mathcal{G}$ is a metric graph, $f$  is a real function and $\D$ is either the Neumann or the Dirichlet Laplacian on $\mathcal{G}$ (see Subsection \ref{lamegra}). Accordingly,  for all $t\in\R_+$ $u(\cdot,t)$ satisfies the {\em Kirchhoff conditon} at any vertex $v\in\mathcal{G}\setminus\pa\mathcal{G}$. In addition, if the boundary $\pa \mathcal G$ is nonempty,
$u$ satisfies either the Neumann boundary condition
\begin{equation}\label{neco}
u_x  = 0 \quad\textrm{in}\,\,\pa\mathcal{G}\times \R_+\,,
\end{equation}
or the Dirichlet boundary condition
\begin{equation}\label{dico}
u = 0 \quad\textrm{in}\,\,\pa\mathcal{G}\times \R_+\,.
\end{equation}

Solutions of problem \eqref{CN} are considered both in the $L^2$ and in the $L^\infty$ framework (see Definitions \ref{deso} and \ref{sol2}; observe that the two solution concepts are different if $\mathcal G$ is infinite).

The following will always be assumed (see Definitions \ref{deg}-\ref{sig}):
$$
\mathcal{G} \;\text{locally finite and connected;}
\leqno(H_0)
$$
$$
 f  \;\;\text{locally Lipschitz continuous}\,.
  \leqno (H_1)
$$
It is not restrictive to suppose $f(0)=0$, as we do in the following.

\smallskip

Monotonicity methods are a major tool for addressing elliptic and parabolic problems in domains of the Euclidean space (in particular, to study attractivity properties of stationary solutions of the latter). It is the purpose of this note to extend such methods to metric graphs satisfying assumption $(H_0)$. The proofs adapt to the present situation ideas which go back to \cite{ACP, dMST}, using a notion of sub- and supersolutions borrowed from the literature concerning impulsive differential problems (e.g., see \cite{LBS, RT1}). In the $L^\infty$ case the proof  relies on the so-called {\it duality method}, already used for similar parabolic equations in the Euclidean space or on Riemannian manifolds (e.g., see  \cite{ACP, GMP, MuP}). Due to the geometry of metric graphs, a nontrivial tailoring of the method to the present situation, as well as  additional structure assumptions on $\mathcal G$ (see $(H_2)$), is needed. We limit ourselves to the simplest semilinear case, yet analogous results could be similarly proven for a wide class of parabolic problems.

Let us mention that, to our knowledge, a direct proof of the maximum principle on metric graphs has only been given in the finite case (see \cite{vb}). The strong maximum principle for the heat equation on infinite graphs, regarded as the irreducibility of the heat semigroup, has been treated in \cite{Ca}.

The mathematical framework is outlined in Section \ref{spatre} and in the Appendix. The main results are stated in Section \ref{prore} and proven in Section \ref{propro}. The case of regular metric trees is discussed in Section \ref{semetre}.


\section{Preliminaries}\label{spatre}
\setcounter{equation}{0}

We refer the reader to the Appendix for notations and  general results used hereafter. Let us only recall that the sets of edges and vertices of $\mathcal{G}$ are denoted by $E$ and $V$, respectively.
\begin{remark}\label{mark}
If $\pa\mathcal{G}= \emptyset$, every statement concerning the set $\mathcal G\setminus\pa\mathcal G$ (respectively $\pa\mathcal G$) holds for the whole of $\mathcal G$ (is void, respectively), and every sum over the elements of (a subdomain of) $\pa\mathcal{G}$ is meant to be zero.
\end{remark}

\subsection{Concepts of solution}\label{spe}

Solutions of problem \eqref{CN} in the $L^2$ framework are meant in the following sense.
\begin{definition}\label{deso}
Let $u_0\in L^2(\mathcal{G})$. By a {\em solution} of problem \eqref{CN} in $[0,T]$ $(T\in\R_+)$ we mean a function $u\in C([0,T];L^2(\mathcal{G}))\cap C^1((0,T];L^2(\mathcal{G}))$, such that $\D u \in C((0,T];L^2(\mathcal{G}))$ and
$$
\left\{
\begin{array}{ll}
u_t(\cdot,t) = (\Delta u)(\cdot,t)+ f(u(\cdot,t))  \quad\text{a.e. in $\mathcal{G}$ for any $t\in (0,T]$\,,}
\\& \\
u(\cdot,0)= u_0\quad\text{a.e. in $\mathcal{G}$}.
\end{array}
\right.
$$
A solution of \eqref{CN} is {\em global}, if it is a solution in $[0,T]$ for any $T\in\R_+$\,.
\end{definition}

\begin{remark}\label{condege}
By Definition \ref{deso}, for any $t\in(0,T]$ $u(\cdot,t)$ belongs to the domain of the Laplacian. Then by \eqref{lave} for all $t\in (0,T]$ there holds $u(\cdot,t)=\bigoplus_{e\in E} u_e(\cdot,t)$, and
\begin{subequations}\label{solge}
\begin{equation}\label{solge0}
 u(\cdot,t)\in H^1(\mathcal{G})\,,\quad  u_e(\cdot,t)\in H^2(I_e) \,\textrm{ for all $e\in E$}, \;
 \quad  \sum_{e\in E}\int_0^{l_e}| u_{exx}(\cdot,t)|^2\,dx<\infty \,,
\end{equation}
\begin{equation}\label{solge1}
u_{et}(\cdot,t)  \,=\,  u_{exx}(\cdot,t) + f(u_e(\cdot,t))  \quad\text{a.e. in $I_e$ for all $e\in E$\,,}
\end{equation}
\begin{equation}\label{solge2}
\sum_{e\ni v}\frac{d u_e}{d\n}(v,t)=0 \quad\text{for all $v\in \mathcal{G}\setminus \pa\mathcal{G}$}\,,
\end{equation}
\begin{equation}\label{solge3}
\text{if $\pa\mathcal{G} \neq \emptyset$, for all $v\in  \pa\mathcal{G}$ }
\quad\left\{
\begin{array}{ll}
u_x(v,t)\,=\,0 \;\;\text{in the Neumann case,}
\\& \\
u(v,t)\,=\,0 \;\;\text{in the Dirichlet case}.
\end{array}
\right.
\end{equation}
\end{subequations}
\end{remark}

\begin{remark}\label{qaz1}
$(i)$ Since $H^1(\mathcal{G})\subseteq C(\mathcal{G})$,  by \eqref{solge0} for all $t\in(0,T]$ $u(\cdot,t)$ is continuous in $\mathcal{G}$, thus in particular at any vertex $v\in \mathcal{G}\setminus  \pa\mathcal{G}$.

\noindent $(ii)$ Since $H^2(I_e)\subseteq C^1(\overline{I}_e)$, there holds $u_e(\cdot,t)\in C^1(\overline{I}_e)$ $(e\in E)$, thus in particular $u_{ex}(\cdot,t)$ is continuous (from one side) at any vertex $v$.

\noindent $(iii)$ Recall that the maps $ i:E\mapsto V$ and $j:\{e\in E\,|\, l_e<\infty\} \mapsto V$ define the initial point and the final point of an edge (see Definition \ref{meg}). If both the inbound star $\Sig_v^+$ and the outbound star $\Sig_v^-$ are nonempty, condition \eqref{solge2} reads:
\begin{equation}\label{tolex1}
\sum_{k=1}^{d_v^+}u_{{e_k}x}(j(e_k),t)\,=\,\sum_{l=1}^{d_v^-} u_{{e_l}x}(i(e_l),t) \quad\text{for all $v\in \mathcal{G}\setminus \pa\mathcal{G}$, $e_k\subseteq\Sig_v^+$, $e_l\subseteq\Sig_v^-$ and $t\in(0,T]$\,.}
\end{equation}
Here $d_{v}^+$ and $d_{v}^-$ are the inbound degree, respectively the outbound degree of the vertex $v$; accordingly, the edges $e_k$ in the left-hand side of \eqref{tolex1} belong to $\Sig_v^+$, the edges $e_l$ in the right-hand side to $\Sig_v^-$ (see Definition \ref{deg}). If either $\Sig_v^+=\emptyset$ or $\Sig_v^-=\emptyset$ $(v\in\mathcal{G})$, the corresponding sum in \eqref{tolex1} is meant to be zero.
\end{remark}

\begin{definition}\label{destage}
A function $q:\mathcal{G}\mapsto \R$, $q\in D(\D)$ is a {\em stationary solution} of  equation \eqref{diffeq} if
\begin{equation}\label{stage}
\D q + f(q) =0  \quad\text{a.e. in $\mathcal{G}$.}
\end{equation}
\end{definition}
\noindent Equivalently, a stationary solution $q=\bigoplus_{e\in E} q_e$ of \eqref{diffeq} satisfies the following:
\begin{subequations}\label{proge}
\begin{equation}\label{proge1}
q\in H^1(\mathcal{G})\,,\quad  q_e\in H^2(I_e) \,\textrm{ for all $e\in E$}, \;
 \quad  \sum_{e\in E}\int_0^{l_e}| q_e''|^2\,dx<\infty \,,
\end{equation}
\begin{equation}\label{proge2}
 q_e''+f(q_e)\,=\,0  \quad\text{a.e. in $I_e$ for all $e\in E$\,,}
\end{equation}
\begin{equation}\label{proge3}
\sum_{e\ni v}\frac{d q_e}{d\n}(v)=0 \quad\text{for all $v\in \mathcal{G}\setminus \pa\mathcal{G}$}\,,
\end{equation}
\begin{equation}\label{proge4}
\text{if $\pa\mathcal{G} \neq \emptyset$, for all $v\in  \pa\mathcal{G}$ }
\quad\left
\{
\begin{array}{ll}
q'(v)\,=\,0 \;\;\text{in the Neumann case,}
\\& \\
q(v)\,=\,0 \;\;\text{in the Dirichlet case}.
\end{array}
\right.
\end{equation}
\end{subequations}

Set $\O_t:=\mathcal{G}\times(0,t)$ $(t\in(0,T])$. In view of Remarks \ref{condege}-\ref{qaz1}, solutions of problem \eqref{CN} in the $L^\infty$ framework are defined as follows:
\begin{definition}\label{sol2}
Let $u_0\in C(\mathcal{G})\cap L^\infty(\mathcal{G)}$.  A function $ u\in C(\overline{\O}_T)\cap L^\infty(\O_T), u(\cdot,t)=\bigoplus_{e\in E} u_e(\cdot,t)$ is said to be a {\em bounded solution}  of problem \eqref{CN} if :

\noindent $(i)$ $u(x, 0)= u_0(x)$ for any $x \in\mathcal G;$

\noindent $(ii)$ for all $e\in E$ there holds
\begin{subequations}\label{mnbv}
\begin{equation}\label{mnbv0}
\text{ $u_e(\cdot, t)\in C^2(I_e)\cap C^1(\bar{I}_e)$ for any $ t\in (0, T]\,,$}
\end{equation}
\begin{equation}\label{mnbv1}
\text{ $u_e(x, \cdot)\in C^1((0, T])$ for any $x\in \bar{I}_e \,;$}
\end{equation}
\end{subequations}

\noindent $(iii)$  equalities \eqref{solge1}-\eqref{solge3} hold.
\end{definition}

Observe that bounded solutions are solutions in the sense of Definition \ref{deso}, if $\mathcal G$ is finite
(similarly for bounded sub- and supersolutions; see Definition \ref{sosol2} below).

\smallskip


\subsection{Subsolutions and supersolutions}\label{wepomo}

Sub- and supersolutions of problem \eqref{CN} in the $L^2$ framework are defined as follows.
\begin{definition}\label{evoss}
A function $\underline{u}\in C([0,T];L^2(\mathcal{G}))\cap\, C^1((0,T];L^2(\mathcal{G}))\cap C((0,T];H^1(\mathcal{G}))$, $\underline{u}(\cdot,t)=\bigoplus_{e\in E} \underline{u}_e(\cdot,t)$ $(t\in (0,T])$, is a {\em subsolution} of problem \eqref{CN} in $[0,T]$ if:
\begin{subequations}\label{sottop}
\begin{equation}\label{sottop0}
 \underline{u}_e(\cdot,t)\in H^2(I_e) \,\textrm{ for all $e\in E$}, \;
 \quad  \sum_{e\in E}\int_0^{l_e}| \underline{u}_{exx}(\cdot,t)|^2\,dx<\infty \,,
\end{equation}
\begin{equation}\label{sottop1}
\underline{u}_{et}(\cdot,t)\,\le\,  \underline{u}_{exx}(\cdot,t)+  f(\underline{u}_e(\cdot,t))  \quad\text{a.e. in $I_e$ for any $e\in E$ and $t\in (0,T]$\,,}
\end{equation}
\begin{equation}\label{solprop22}
\underline{u}(\cdot,0)\,\le\, u_0 \;\;\text{a.e. in $\mathcal{G}$}\,,
\end{equation}
\begin{equation}\label{solprop23}
\sum_{e\ni v}\frac{d \underline{u}_e}{d\n}(v,t)\;\le\;0 \quad\text{for all $v\in \mathcal{G}\setminus \pa\mathcal{G}$ and $t\in(0,T]$}\,,
\end{equation}
\begin{equation}\label{solprop3}
\text{if $\pa\mathcal{G} \neq \emptyset$, for all $(v,t)\in  \pa\mathcal{G}\times(0,T]$ }
\quad
\left\{
\begin{array}{ll}
&\underline{u}_{ex}(v,t)\le0 \;\text{if $v=j(e)$} , \;\; \underline{u}_{ex}(v,t)\ge0 \;\text{if $v=i(e)$}\\
&\text{in the Neumann case,} \smallskip
\\ 
&u(v,t)\,\le\,0 \;\;\text{in the Dirichlet case}.
\end{array}
\right.
\end{equation}
\end{subequations}
A {\em supersolution} $\overline{u}$ of \eqref{CN} in $[0,T]$ is defined by reversing inequalities in \eqref{sottop1}-\eqref{solprop3}.
\end{definition}

\begin{remark}\label{limiss}
$(i)$ Since $H^1(\mathcal{G})\subseteq C(\mathcal{G})$ and $H^2(I_e)\subseteq C^1(\overline{I}_e)$, by \eqref{sottop0} for all $t\in(0,T]$:

\noindent $(a)$ both $\underline{u}(\cdot,t)$ and $\overline{u}(\cdot,t)$ are continuous in $\mathcal{G}$, thus at any vertex $v\in \mathcal{G}\setminus  \pa\mathcal{G}$, and there holds $\underline{u}(\cdot,t)\in L^\infty_{loc}(\mathcal{G})$, $\overline{u}(\cdot,t)\in L^\infty_{loc}(\mathcal{G})$;

\noindent $(b)$ there holds $\underline{u}_e(\cdot,t), \overline{u}_e(\cdot,t)\in C^1(\overline{I}_e)$, thus in particular $\underline{u}_{ex}(\cdot,t)$ and $\overline{u}_{ex}(\cdot,t)$ are continuous (from one side) at any vertex $v$.

\noindent $(ii)$ Condition \eqref{solprop23} can be rewritten as follows:
\begin{equation}\label{solex1}
\sum_{k=1}^{d_v^+} \underline{u}_{{e_k}x}(j(e_k),t)\,\le\,\sum_{l=1}^{d_v^-} \underline{u}_{{e_l}x}(i(e_l),t) \quad\text{for all $v\in \mathcal{G}\setminus \pa\mathcal{G}$, $e_k\subseteq\Sig_v^+$, $e_l\subseteq\Sig_v^-$ and $t\in(0,T]$}\,.
\end{equation}
\end{remark}

\begin{definition}\label{soso}
A function $\underline{q}:\mathcal{G}\mapsto \R$, $\underline{q}=\bigoplus_{e\in E} \underline{q}_e,$ is a {\em stationary subsolution} of equation \eqref{diffeq} if
\begin{subequations}\label{sproge}
\begin{equation}\label{sproge1}
\underline{q}\in H^1(\mathcal{G})\,,\quad  \underline{q}_e\in H^2(I_e) \,\textrm{ for all $e\in E$}, \;
 \quad  \sum_{e\in E}\int_0^{l_e}| \underline{q}_e''|^2\,dx<\infty \,,
\end{equation}
\begin{equation}\label{sproge2}
 \underline{q}_e''+f(\underline{q}_e)\,\ge\,0  \quad\text{a.e. in $I_e$ for all $e\in E$\,,}
\end{equation}
\begin{equation}\label{sproge3}
\sum_{e\ni v}\frac{d  \underline{q}_e}{d\n}(v)\;\le\;0 \quad\text{for all $v\in \mathcal{G}\setminus \pa\mathcal{G}$}\,,
\end{equation}
\begin{equation}\label{sproge4}
\text{if $\pa\mathcal{G} \neq \emptyset$, for all $v\in  \pa\mathcal{G}$}
\quad\left\{
\begin{array}{ll}
 \underline{q}_e'(v)\le0 \;\;\text{if $v=j(e)$} , \;\; \underline{q}_e'(v)\ge0\;\; \text{if $v=i(e)$ in the Neumann case,}
\\& \\
 \underline{q}(v)\,\le\,0 \;\;\text{in the Dirichlet case}.
\end{array}
\right.
\end{equation}
\end{subequations}
A {\em stationary supersolution} $\overline{q}$ is defined by reversing inequalities in \eqref{sproge2}-\eqref{sproge4}.
\end{definition}
In agreement with Definition \ref{sol2}, we give a different concept of sub- and supersolutions for the $L^\infty$ framework:
\begin{definition}\label{sosol2}
A function $ u\in C(\overline{\O}_T)\cap L^\infty(\O_T), u(\cdot,t)=\bigoplus_{e\in E} u_e(\cdot,t)$ is said to be a {\em bounded subsolution}  of problem \eqref{CN}, if for all $e\in E$ $u_e$ has the regularity stated in \eqref{mnbv}, and
inequalities \eqref{sottop1}-\eqref{solprop3} hold.
\noindent A {\em bounded supersolution} of \eqref{CN} is similarly defined, with reverse inequalities in \eqref{sottop1}-\eqref{solprop3}.
\end{definition}

By abuse of language, when we say ``sub-, supersolutions or solutions to problem \eqref{CN}'' we always mean sub-, supersolutions or solutions in the $L^2$ framework (see Definitions \ref{deso} and \ref{evoss}). Sub-, supersolutions or solutions in the $L^\infty$ framework, on the other hand, are always explicitly referred to as ``bounded sub-, supersolutions or solutions'' (see Definitions \ref{sol2} and \ref{sosol2}).


\section{Main results}\label{prore}
\setcounter{equation}{0}

The following result concerning well-posedness  of problem \eqref{CN} (as well as a companion result concerning bounded solutions, whose formulation is left to the reader) follows by standard arguments of semigroup theory (e.g., see \cite[Theorem 4.4]{Y}).
\begin{theorem}\label{wpp}
Let  $u_0\in L^2(\mathcal G)$, and let $(H_0)$-$(H_1)$ hold. Then there exists $T=T(u_0)\in\R_+$ such that problem \eqref{CN} has a unique solution in $[0,T]$. In addition,
 if $u_0\ge0$ in $\mathcal{G}$ and $f:\overline{\R}_+\mapsto\overline{\R}_+$,  there holds $u(\cdot,t)\ge0$ in $\mathcal{G}$  for any $t\in[0,T]$.
\end{theorem}

The following comparison result will be proven.
\begin{theorem}\label{geco}
Let $(H_0)$-$(H_1)$ hold. Let $\underline{u}, \,\overline{u}$ be a subsolution, respectively a supersolution of problem \eqref{CN} in $[0,T]$. Then there holds $\underline{u} \le \overline{u}$ in $\O_T$\,.
\end{theorem}

As a consequence of Theorem \ref{geco}, the following holds:
\begin{theorem}\label{decre}
Let $(H_0)$-$(H_1)$ hold. Let $\underline{q},\,\overline{q}$ be a stationary sub-, respectively supersolution of  \eqref{diffeq}, such that $\underline{q}\le\overline{q}$ a.e. in $\mathcal{G}$. Let $u_1$, $u_2$ be the solutions of \eqref{CN} with $u_0=\underline{q},$ respectively $u_0=\overline{q}$. Then:

\noindent $(i)$ $u_1$, $u_2$ are global, and for all $t\in\R_+$ there holds $\underline{q}\le u_1(\cdot,t)\le u_2(\cdot,t)\le\overline{q}$ in $\mathcal{G}$;

\noindent $(ii)$ if $\underline{q}\le u_0\le\overline{q}$, the solution $u$ of problem \eqref{CN}  is global, and for all $t\in\R_+$ there holds $u_1(\cdot,t)\le u(\cdot,t)\le u_2(\cdot,t)$ in $\mathcal{G}$;

\noindent $(iii)$ the map $t\mapsto u_1(\cdot,t)$ is nondecreasing, $t\mapsto u_2(\cdot,t)$ is nonincreasing, and the pointwise limits $\hat{u}_1:=\lim_{t\to\infty}u_1(\cdot,t)$, $\hat{u}_2:=\lim_{t\to\infty}u_2(\cdot,t)$ are stationary solutions of \eqref{diffeq}. Moreover, there holds $\hat{u}_1\le \hat{u}_2$\,, and every stationary solution of \eqref{diffeq} with $\underline{q}\le \tilde{u}\le\overline{q}$ satisfies $\hat{u}_1\le\tilde{u}\le\hat{u}_2$\,.
\end{theorem}
The last statement in $(ii)$ above is referred to by saying that $\hat{u}_1,\, \hat{u}_2$ are the {\em minimal}, respectively the {\em maximal} stationary solution of \eqref{diffeq} in the function interval $[\underline{q},\overline{q}]:=\{q\in L^2(\mathcal{G})\,|\,\underline{q}\le q\le\overline{q}\;\text{in $\mathcal{G}$}\}$.
\begin{remark}\label{dinico}
By the regularity of $u_i(\cdot,t)$ and $\hat{u}_i$, and the monotone character of the map $t\mapsto u_i(\cdot,t)$, it follows easily from  the Dini Theorem that the convergence $u_i(\cdot,t)\to\hat{u}_i$ in Theorem \ref{decre} is uniform in any  compact subset of $\mathcal{G}$ $(i=1,2)$.
\end{remark}
\smallskip

Theorem \ref{comparison} below is the counterpart of Theorem \ref{geco} for bounded sub- and supersolutions; its statement needs some preliminary remarks.
Let $v_0\in V$ be arbitrarily fixed, and let $r(x):=d(x, v_0)$ denote the distance between $v_0$ and any $x\in \mathcal G$ (see Subsection \ref{megra}).
For any $R>0$ set
$$
V_R:=\{v\in V\,|\, r(v) \leq R \}\,, \quad S_R:= \{v\in V\,|\, r(v)= R\}\,,\quad
$$
$$
E_R:=\{e\in E\,|\, i(e)\in V_R, j(e)\in V_R\}\,, \quad \mathcal G_R:=(V_R, E_R)\subseteq \mathcal G\,.
$$
For any sequence $\{R_n\}\subset \R_+$ we denote $V_n\equiv V_{R_n}$, and similarly for
$S_n,\,E_n,\, \mathcal G_n$ $(n\in\N)$.
\smallskip

When proving Theorem \ref{comparison}, we choose an orientation on $\mathcal G$ such that for any $e\equiv(v_1, v_2)\in E$ with $r(v_1)\leq r(v_2)$ there holds $v_1=i(e)$ and $v_2=j(e).$ Moreover, we assume that
the {\it jump size} of $\mathcal G$ is finite, namely
$$
s\,:=\,\sup\{d(i(e), j(e))<\infty\,|\, e\in E\}<\infty\,.
$$
Therefore, we can fix a diverging sequence $\{R_n\}\subset \R_+$ such that
\begin{equation}\label{e30f}
\mathcal G_{m}\subset \mathcal G_{n} \;\, \text{ if } m<n, \quad
\bigcup_{n=1}^{\infty}\mathcal G_{n}=\mathcal G,
\frac 1{c_0}\leq R_n-R_{n-1}\leq c_0\;\, \text{ for  some $c_0>1$ }\qquad (n\in \mathbb N)\,.
\end{equation}

For any $v\in V_n$ we set $\Sig_{v,n}^\pm:= \Sig_v^\pm\cap\mathcal G_n$, with $\Sig_v^\pm$ defined by \eqref{sv1}. Accordingly, we denote by $d_{v,n}^+$ and $d_{v,n}^-$ the inbound, respectively outbound degree of $v$ relative to $\mathcal G_n$ - namely, the number of edges $e\in E_n$ with $j(e)=v$ (respectively, with $i(e)=v$). Due to the orientation chosen on $\mathcal G$ and the definition of $\mathcal G_n$\,, if $v$ belongs to $V_n\setminus S_n$ every edge $e$ with either $j(e)=v$ or $i(e)=v$ belongs to $E_n$\,, hence there holds $d_{v,n}^\pm= d_v^\pm$\,. It follows that:
\begin{subequations}\label{e1f}
\begin{eqnarray}\label{e1f1}
\sum_{e\in E_n, e\ni v}\frac{d f_e}{d\n}(v)&=& \sum_{k=1}^{d_{v,n}^+}f_{e_k}'(j(e_k)) \;- \sum_{l=1}^{d_{v,n}^-}f_{e_l}'(i(e_l))
\,=\\
&=& \sum_{k=1}^{d_v^+}f_{e_k}'(j(e_k)) \;- \sum_{l=1}^{d_v^-}f_{e_l}'(i(e_l)) \quad \text{ for any } v\in V_n\setminus S_n\,. \nonumber
\end{eqnarray}
Similarly, every $v\in S_n$ is the final point of some $e\in E_n$, thus there holds $d_{v,n}^+= d_v^+$, $d_{v,n}^-=0.$ Therefore,
\begin{equation}\label{e1f2}
\sum_{e\in E_n, e\ni v}\frac{d f_e}{d\n}(v)\,=\, \sum_{k=1}^{d_{v,n}^+}f_{e_k}'(j(e_k))\,=\, \sum_{k=1}^{d_v^+}f_{e_k}'(j(e_k)) \quad \text{ for any } v\in S_n\,.
\end{equation}
\end{subequations}

We can now state the following assumption:
\begin{equation*}
\left\{\begin{array}{ll}
(i) & d^+_v\,\le\, d^-_v \; \text{ for any } v\in V;\smallskip\\
(ii)&\text{there exist $ C>0,\, \theta>0,\, \beta\in [0, 2]$ such that } \smallskip\\ & \sum_{v\in S_n}d_v^+ \leq  C \exp\{\theta   R_n^{\beta}\}\;\; \text{for any} \;\; n\in\N\,.
\end{array}
\right.
\leqno(H_2)
\end{equation*}
\begin{remark}
It is easily seen that assumption $(H_2)$ is satisfied for any regular tree $\mathcal{T}$ (see Subsection \ref{metre}). In fact,
in this case the boundary $\pa\mathcal{T}$ only consists of the root $O$, and a natural choice is to assume $v_0 =O$. For any vertex $v$ of the $n$-th generation (thus belonging to $\mathcal{T}\setminus O$) there holds $d_v^+=1$, $d_v^-=b_n\ge2$, hence $(H_2)$-$(i)$ is satisfied. Concerning $(H_2)$-$(ii)$, assume for simplicity that $\mathcal T$ is homogeneous. Then for any vertex $v$ of the $n$-th generation there holds $b_n=b$, $r(v)=nr$, thus choosing $R_n=nr$ gives $\sum_{v\in S_n}d_v^+=b^{n-1}$. It is easily seen that $(H_2)$-$(ii)$ holds in this case with $C\ge\frac1b$, $\theta\ge\frac{\log{b}}{r}$ and $\b=1$.

\end{remark}
\begin{theorem}\label{comparison}
Let $\mathcal G$ be infinite with finite jump size, and let $(H_0)$-$(H_2)$ hold. Let one of the following assumptions be satisfied:

\noindent $(i)$ $\pa\mathcal G=\emptyset\,;$

\noindent $(ii)$ $\pa\mathcal G\neq\emptyset$ is finite and the Neumann boundary condition \eqref{neco} holds;

\noindent $(iii)$ $\pa\mathcal G\neq\emptyset$ and the Dirichlet boundary condition \eqref{dico} holds.

\noindent Let $\underline u$, $\overline u$ be a bounded subsolution, respectively a bounded supersolution of \eqref{CN}. Then there holds  $\underline{u} \le \overline{u}$ in $\O_T$\,.
\end{theorem}

\begin{remark}\label{pagifi}
Observe that assumption $(ii)$ of Theorem \ref{comparison} holds for regular trees with Neumann boundary condition at the root $O$ (see Section \ref{semetre}); in this case a natural choice is $v_0=O$\,.
\end{remark}

A companion result of Theorem \ref{decre} follows from Theorem \ref{comparison}; we leave its formulation to the reader.


\section{Proofs}\label{propro}
\setcounter{equation}{0}

\noindent {\em Proof of Theorem \ref{geco}.} $(i)$ Set $w:=\underline{u}-\overline{u}$, and for all $(x,t)\in \O_T$
$$
 \qquad a=a(x,t):=
\left\{
\begin{array}{ll}
\frac{f(\underline{u}(x,t))-f(\overline{u}(x,t))}{w(x,t)} \quad\text{if $w(x,t)\neq0$\,,}
\\& \\
0 \quad \text{otherwise}.
\end{array}
\right.
$$
Observe that $a\in L^\infty_{loc}(\O_T)$. In fact,  by assumption $(H_1)$ there holds
\begin{equation}\label{LM}
|a(x,t)|\le L_M\quad\text{for any compact $K\subseteq \O_T$}\qquad ((x,t)\in K)\,,
\end{equation}
 $L_M$ denoting the Lipschitz constant for $f$ in the interval $[-M,M]$ with $M\equiv M(K):=\max\{\|\underline{u}\|_{L^\infty(K)}\,,\|\overline{u}\|_{L^\infty(K)}\}$ (observe that both $\underline{u}$ and $\overline{u}$ belong to $L^\infty_{loc}(\O_T)$, since $\underline{u},\overline{u}\in C((0,T];H^1(\mathcal{G}))\subseteq C((0,T];C(\mathcal{G}))$).

By Definition \ref{evoss} there holds
\begin{subequations}\label{fine}
\begin{equation}\label{fine11}
w(\cdot,0)\,\le\, 0 \;\;\text{a.e. in $\mathcal{G}$}\,,
\end{equation}
and  for any $t\in(0,T]$ (see Remark \ref{limiss}-$(ii)$):
\begin{equation}\label{fine0}
w(\cdot,t)\in H^1(\mathcal{G})\,,\quad w_e(\cdot,t)\in H^2(I_e) \,\textrm{ for all $e\in E$}, \;
 \quad  \sum_{e\in E}\int_0^{l_e}| w_{exx}(\cdot,t)|^2\,dx<\infty \,,
\end{equation}
\begin{equation}\label{fine1}
w_t(\cdot,t)-w_{xx}(\cdot,t) - a(\cdot,t)w(\cdot,t)\,\le\, 0 \quad\text{a.e. in $\mathcal{G}$} \,,
\end{equation}
\begin{equation}\label{fine2}
\sum_{k=1}^{d_v^+} w_{{e_k}x}(v,t)\,\le\,\sum_{l=1}^{d_v^-} w_{{e_l}x}(v,t) \quad\text{for all $v\in  \mathcal{G}\setminus\pa \mathcal{G}$, $e_k\subseteq\Sig_v^+$, $e_l\subseteq\Sig_v^-$}\,,
\end{equation}
\begin{equation}\label{fine21}
\text{if $\pa\mathcal{G}= \emptyset$,  for all $(v,t)\in  \pa\mathcal{G}\times(0,T]$ }
\quad\left\{
\begin{array}{ll}
&w_{ex}(v,t)\le0 \;\;\text{if $v=j(e)$} , \;\; w_{ex}(v,t)\ge0\;\; \text{if $v=i(e)$} \\
& \text{in the Neumann case,}\smallskip
 \\
&w(v,t)\,\le\,0 \;\;\text{in the Dirichlet case}.
\end{array}
\right.
\end{equation}
\end{subequations}
Observe that by \eqref{fine0} for all $t\in(0,T]$ there holds $w(\cdot,t)\in C(\mathcal{G})$, thus in particular $w(\cdot,t)$ is continuous at each vertex $v\in \mathcal{G}\setminus\pa\mathcal{G}$.

Let $\t\in(0,T]$ be fixed. Consider the backward problem
\begin{equation*}
\left\{
\begin{array}{ll}
\varphi_t \,=\, -\D \varphi  -a\varphi \quad\textrm{in}\,\,\O_\t
\\& \\
\varphi = \z
\quad \textrm{in}\,\, \mathcal{G}\times \{\t\} \,,
\end{array}
\right.
\leqno{(BP)}
\end{equation*}
where $\z=\bigoplus_{e\in E} \z_e\in C_0^\infty(\mathcal{G})$, $0\le \z\le 1$. Since the Laplacian generates an analytic semigroup in $L^2(\mathcal{G})$, by standard results there exists a unique nonnegative solution $\varphi=\bigoplus_{e\in E} \varphi_e\in C([0, \t];L^2( \mathcal{G}))\cap C^1([0, \t);L^2( \mathcal{G}))$, with $\D \varphi \in C([0, \t);L^2( \mathcal{G}))$, of problem $(BP)$.
In particular, there holds $\varphi (\cdot,t)\in L^2(\mathcal{G})$, $\varphi (\cdot,t)\ge0$ for any $t\in[0,\t]$, and  for all $t\in(0,\t):$
\begin{subequations}
\begin{equation}\label{duffop1}
\varphi(\cdot,t)\in H^1(\mathcal{G})\,, \quad \varphi_e(\cdot,t)\in H^2(I_e) \,\textrm{ for all $e\in E$}\,,
\end{equation}
\begin{equation}\label{duffop2}
\varphi_t(\cdot,t) +\varphi_{xx}(\cdot,t)+a(\cdot,t)\varphi(\cdot,t)=0
\quad\text{a.e. in $\mathcal{G}$\,,}
\end{equation}
\begin{equation}\label{duffop3}
\sum_{k=1}^{d_v^+}\varphi_{{e_k}x}(j(e_k),t)\,=\,\sum_{l=1}^{d_v^-} \varphi_{{e_l}x}(i(e_l),t) \quad\text{for all $v\in \mathcal{G}\setminus \pa\mathcal{G}$, $e_k\subseteq\Sig_v^+$, $e_l\subseteq\Sig_v^-$\,,}
\end{equation}
\begin{equation}\label{duffop4}
\text{if $\pa\mathcal{G}= \emptyset$, for all $v\in  \pa\mathcal{G}$ }
\quad\left\{
\begin{array}{ll}
\varphi_x(v,t)\,=\,0 \;\;\text{in the Neumann case,}
\\& \\
\varphi(v,t)\,=\,0 \;\;\text{in the Dirichlet case}.
\end{array}
\right.
\end{equation}
\end{subequations}
In particular, by \eqref{duffop1} for all $t\in[0,\t)$ there holds $\varphi(\cdot,t)\in C(\mathcal{G})$, thus in particular $\varphi(\cdot,t)$ is continuous at each vertex $v\in \mathcal{G}\setminus\pa\mathcal{G}$.

\smallskip

Since $\varphi\ge0$, from inequality \eqref{fine1} for all $e\in E$ we get
\begin{equation}\label{ar1}
 \int_0^\t dt\int_0^{l_e}dx\left[w_{et}(x,t)-w_{exx}(x,t) - a_e(x,t)w_e(x,t)\right]\varphi_e(x,t)\,\le\,0\,,
 \end{equation}
where
\begin{equation}\label{arie}
 \qquad a_e(x,t):=
\left\{
\begin{array}{ll}
\frac{f(\underline{u}_e(x,t))-f(\overline{u}_e(x,t))}{w_e(x,t)} \quad\text{if $w_e(x,t)\neq0$\,,}
\\& \\
0 \quad \text{otherwise}.
\end{array}
\right.
 \end{equation}

Below we shall prove that for any $e\in E$ there holds
\begin{eqnarray}\label{claim}
&&
 \int_0^\t dt\int_0^{l_e}dx\left[w_{et}(x,t)-w_{exx}(x,t) - a_e(x,t)w_e(x,t)\right]\varphi_e(x,t)\,\ge\\
&\ge& \int_0^{l_e}w_e(x,\t)\,\z_e(x)\,dx \,-\nonumber
\\
&-& \int_0^\t \left[ w_{ex}(j(e),t)\varphi_{e}(j(e),t)\,-\, w_{ex}(i(e),t)\varphi_{e}(i(e),t)\right]\,dt  \,+\nonumber\\
&+&  \int_0^\t \left[ w_{e}(j(e),t)\varphi_{ex}(j(e),t) \,-\,w_{e}(i(e),t)\varphi_{ex}(i(e),t)
\right]\,dt  \,.\nonumber
\end{eqnarray}

From \eqref{ar1} and \eqref{claim} we obtain
\begin{eqnarray}\label{pos1}
0&\ge&
\int_0^\t dt\int_{\mathcal{G}}dx\left[w_t(x,t)-w_{xx}(x,t) - a(x,t)w(x,t)\right]\varphi(x,t)\,\ge\\
&\ge& \int_{\mathcal{G}}w(x,\t)\,\z(x)\,dx\,-\ \nonumber
\\
&-&\sum_{e\in E}
\int_0^\t \left[ w_{ex}(j(e),t)\varphi_{e}(j(e),t)\,-\, w_{ex}(i(e),t)\varphi_{e}(i(e),t)
\right]\,dt  \,+\nonumber\\
&+&  \sum_{e\in E}
\int_0^\t
\left[ w_{e}(j(e),t)\varphi_{ex}(j(e),t) \,-\,w_{e}(i(e),t)\varphi_{ex}(i(e),t)
\right]\,dt  \,.\nonumber
\end{eqnarray}

If $\pa\mathcal{G}\neq \emptyset$, recall that $d_v=1$ for all $v\in \pa \mathcal{G}$ and set
$$
\pa \mathcal{G}^\pm:= \{v\in V\,|\, d_v^\pm=1\,, \; d_v^\mp=0\}\subseteq\pa \mathcal{G}\quad \Rightarrow\quad
\pa \mathcal{G}=\pa \mathcal{G}^+\cup\pa \mathcal{G}^-\,, \;\; \pa \mathcal{G}^+\cap\pa \mathcal{G}^-=\emptyset\,.
$$
Also observe that for any $t\in[0,\t)$, $v\in \mathcal{G}\setminus\pa\mathcal{G}$ and $e\in E$ with $e\ni v$ there holds $\varphi_{e}(v,t)=\varphi(v,t)$, since $\varphi(\cdot,t)$ is continuous at any $v\in \mathcal{G}\setminus\pa\mathcal{G}$. Then we get plainly:
\begin{eqnarray*}
&&\sum_{e\in E}
\int_0^\t \left[ w_{ex}(j(e),t)\varphi_{e}(j(e),t)\,-\, w_{ex}(i(e),t)\varphi_{e}(i(e),t)
\right]\,dt  \,=\\
&=&  \sum_{v\in \mathcal{G}\setminus\pa\mathcal{G}}
\int_0^\t
\Big[ \sum_{k=1}^{d_v^+}w_{{e_k}x}(j(e_k),t)\varphi_{e_k}(j(e_k),t) \,-\, \sum_{l=1}^{d_v^-}w_{{e_l}x}(i(e_l),t)\varphi_{e_l}(i(e_l),t)
\Big]\,dt  \,+\nonumber\\
&+&\int_0^\t \Big[\sum_{v\in \pa\mathcal{G}^+} w_x(v,t)\varphi(v,t)\,-\,\sum_{v\in \pa\mathcal{G}^-} w_x(v,t)\varphi(v,t)\Big]\,dt\,=\nonumber\\
&=&  \sum_{v\in \mathcal{G}\setminus\pa\mathcal{G}}
\int_0^\t
\Big[ \sum_{k=1}^{d_v^+}w_{{e_k}x}(v,t) \,-\, \sum_{l=1}^{d_v^-}w_{{e_l}x}(v,t)\Big]\varphi(v,t)\,dt  \,+\nonumber\\
&+&\int_0^\t \Big[\sum_{v\in \pa\mathcal{G}^+} w_x(v,t)\varphi(v,t)\,-\,\sum_{v\in \pa\mathcal{G}^-} w_x(v,t)\varphi(v,t)\Big]\,dt\,,\nonumber
\end{eqnarray*}
whence by \eqref{fine2}
\begin{eqnarray}\label{pos11}
&&\sum_{e\in E}
\int_0^\t \left[ w_{ex}(j(e),t)\varphi_{e}(j(e),t)\,-\, w_{ex}(i(e),t)\varphi_{e}(i(e),t)
\right]\,dt  \,\le\\
&\le&\int_0^\t \Big[\sum_{v\in \pa\mathcal{G}^+} w_x(v,t)\varphi(v,t)\,-\,\sum_{v\in \pa\mathcal{G}^-} w_x(v,t)\varphi(v,t)\Big]\,dt\,.\nonumber
\end{eqnarray}
If $\pa\mathcal{G}\neq \emptyset$, observe that
\begin{equation}\label{pos2}
\int_0^\t \Big[\sum_{v\in \pa\mathcal{G}^+} w_x(v,t)\varphi(v,t)\,-\,\sum_{v\in \pa\mathcal{G}^-} w_x(v,t)\varphi(v,t)\Big]\,dt\,\le\,0\,,
\end{equation}
since in the Neumann case for any $t\in(0,\t)$ there holds $\pm w_x(\cdot,t)\le0$ at $\pa \mathcal{G}^\pm$ (see \eqref{fine21}), whereas in the Dirichlet case $\varphi(\cdot,t)=0$ at $\pa \mathcal{G}$. Therefore, by \eqref{pos11} if $\pa\mathcal{G}= \emptyset$, or by \eqref{pos11}-\eqref{pos2} otherwise, we get
\begin{equation}\label{pos3}
\sum_{e\in E}\int_0^\t \left[ w_{ex}(j(e),t)\varphi_{e}(j(e),t)\,-\, w_{ex}(i(e),t)\varphi_{e}(i(e),t)\right]\,dt  \,\le\,0\,.
\end{equation}

It is similarly seen that
\begin{equation}\label{pos4}
\sum_{e\in E} \int_0^\t \left[ w_{e}(j(e),t)\varphi_{ex}(j(e),t) \,-\,w_{e}(i(e),t)\varphi_{ex}(i(e),t) \right]\,dt \,\ge\,0\,.
\end{equation}
In fact, recalling that for all $t\in[0,\t)$ $w(\cdot,t)$ is continuous at each $v\in \mathcal{G}\setminus\pa\mathcal{G}$,
thus $w_{e}(v,t)=w(v,t)$ for any $e\in E$ with $e\ni v$, and using \eqref{duffop3} we obtain
\begin{eqnarray*}
&&\sum_{e\in E}
\int_0^\t \left[ w_{e}(j(e),t)\varphi_{ex}(j(e),t) \,-\,w_{e}(i(e),t)\varphi_{ex}(i(e),t) \right]\,dt \,=\\
&=&  \sum_{v\in \mathcal{G}\setminus\pa\mathcal{G}}
\int_0^\t
\Big[ \sum_{k=1}^{d_v^+}w_{e_k}(j(e_k),t)\varphi_{{e_k}x}(j(e_k),t) \,-\, \sum_{l=1}^{d_v^-}w_{e_l}(i(e_l),t)\varphi_{{e_l}x}(i(e_l),t)
\Big]\,dt  \,+\nonumber\\
&+&\int_0^\t \Big[\sum_{v\in \pa\mathcal{G}^+} w(v,t)\varphi_x(v,t)\,-\,\sum_{v\in \pa\mathcal{G}^-} w(v,t)\varphi_x(v,t)\Big]\,dt\,=\nonumber\\
&=&  \sum_{v\in \mathcal{G}\setminus\pa\mathcal{G}}
\int_0^\t
\Big[ \sum_{k=1}^{d_v^+}\varphi_{{e_k}x}(j(e_k),t) \,-\, \sum_{l=1}^{d_v^-}\varphi_{{e_l}x}(i(e_l),t) \Big]w(v,t)\,dt  \,+\nonumber\\
&+&\int_0^\t \Big[\sum_{v\in \pa\mathcal{G}^+} w(v,t)\varphi_x(v,t)\,-\,\sum_{v\in \pa\mathcal{G}^-} w(v,t)\varphi_x(v,t)\Big]\,dt\,=\\
&=&\int_0^\t \Big[\sum_{v\in \pa\mathcal{G}^+} w(v,t)\varphi_x(v,t)\,-\,\sum_{v\in \pa\mathcal{G}^-} w(v,t)\varphi_x(v,t)\Big]\,dt\,\ge\,0\,;\nonumber
\end{eqnarray*}
in fact,if $\pa\mathcal G\neq\emptyset$,  in the Dirichlet case there holds $w(\cdot,t)\le0$ at $\pa \mathcal{G}$ (see \eqref{solprop3}) and $\pm \varphi_x(\cdot,t)\le0$ at $\pa \mathcal{G}^\pm$, whereas in the Neumann case by \eqref{duffop4} there holds $\varphi_x(\cdot,t)=0$ at $\pa \mathcal{G}$ $(t\in(0,\t))$.

\smallskip

In view of \eqref{pos1}, \eqref{pos3} and \eqref{pos4}, there holds
\begin{equation}\label{claim3}
 \int_{\mathcal{G}}w(x,\t)\,\z(x)\,dx\,\le\,0 \,.
\end{equation}
By the arbitrariness of $\z$, in the above inequality  we can choose
$$
\z=\z_k\,, \qquad \lim_{k\to\infty}\z_k(x)= \chi_{\{w(\cdot,\t)>0\}}(x) \quad\text{for any $x\in\mathcal{G}$\,.}
$$
Then letting $k\to\infty$ in \eqref{claim3} (written with $\z=\z_k$) gives
\begin{equation*}
\int_{\mathcal{G}}\;[w(x,\t)]_+\,dx \,\le\,0\,,
\end{equation*}
whence by the arbitrariness of $\t\in\R_+$ and the regularity of $w$ the conclusion follows.

\smallskip

It remains to prove \eqref{claim}.
To this purpose,  observe that by  inequality \eqref{fine11}
\begin{eqnarray}\label{term1}
&&\int_0^\t dt\int_0^{l_e}dx\,w_{et}(x,t)\varphi_e(x,t)\,=\\
&=& \int_0^{l_e}w_e(x,\t) \varphi_e(x,\t)\,dx \,-\, \int_0^{l_e}w_e(x,0)\varphi_e(x,0)\,dx \,-\, \int_0^\t dt\int_0^{l_e}dx\,w_e(x,t)\varphi_{et}(x,t)
\,\ge \nonumber\\
&\ge& \int_0^{l_e}w_e(x,\t)\z_e(x)\,dx \,-\, \int_0^\t dt\int_0^{l_e}dx\,w_e(x,t)\varphi_{et}(x,t)\qquad\qquad (e\in E)\,. \nonumber
\end{eqnarray}
On the other hand, for all $e\in E$
\begin{eqnarray}\label{term1110}
&&\int_0^\t dt\int_0^{l_e}dx\,w_{exx}(x,t)\varphi_{e}(x,t)\,= \int_0^\t dt\int_0^{l_e}dx\,w_e(x,t)\varphi_{exx}(x,t)\,+\\
&+& \int_0^\t \big[w_{ex}(l_e,t) \varphi_{e}(l_e,t)- w_{ex}(0,t) \varphi_{e}(0,t)\big]
\,dt\,- \nonumber
\\
&-& \int_0^\t \big[w_e(l_e,t) \varphi_{ex}(l_e,t)- w_e(0,t) \varphi_{ex}(0,t)\big]\,dt\,. \nonumber
\end{eqnarray}
Recalling that $\pi_e^{-1} (0)=i(e)$,  $\pi_e^{-1} (l_e)=j(e)$ (where $\pi_e: \mathcal{G}_e\mapsto I_E$ is the canonical injection; see Subsection \ref{megra}), from \eqref{term1}-\eqref{term1110} we get
\begin{eqnarray*}
&&
 \int_0^\t dt\int_0^{l_e}dx\left[w_{et}(x,t)-w_{exx}(x,t) - a_e(x,t)w_e(x,t)\right]\varphi_e(x,t)\,\ge\\
&\ge& \int_0^{l_e}w_e(x,\t)\,\z_e(x)\,dx\,-\,\int_0^\t dt\int_0^{l_e}dx\left[\varphi_{et}(x,t)+\varphi_{exx}(x,t) + a_e(x,t)\varphi_e(x,t)\right]w_e(x,t) \,-
\nonumber
\\
&-& \int_0^\t \left[ w_{ex}(j(e),t)\varphi_{e}(j(e),t)\,-\, w_{ex}(i(e),t)\varphi_{e}(i(e),t)\right]\,dt  \,+\nonumber\\
&+&  \int_0^\t \left[ w_{e}(j(e),t)\varphi_{ex}(j(e),t) \,-\,w_{e}(i(e),t)\varphi_{ex}(i(e),t)
\right]\,dt  \,,\nonumber
\end{eqnarray*}
which by \eqref{duffop2} coincides with \eqref{claim}. This completes the proof.
\hfill$\square$

\smallskip

\noindent {\em Proof of Theorem \ref{decre}.}
Since $\overline{q}$ is a (stationary) supersolution and $u_2$ is a solution of problem \eqref{CN} with $u_0=\overline{q}$, by Theorem \ref{geco}
there holds $u_2(x,t)\le\overline{q}(x)$ for any $x\in\mathcal{G}$ and $t$ in the maximal interval of existence of $u_2$. Then by a standard continuation arguments $u_2$ is global, and $u_2(x,t)\le\overline{q}(x)$ for any $x\in\mathcal{G}$ and $t\in\R_+$\,. It is similarly seen that  $u_1$ is global, and there holds $u_1(x,t)\ge\underline{q}(x)$ for any $x\in\mathcal{G}$ and $t\in\R_+$\,. The other statements in claims $(i)$-$(ii)$ follow from Theorem \ref{geco}.

Concerning $(iii)$, consider problem \eqref{CN} with initial data $u_2(x,t)$, $t\in\R_+$ arbitarily fixed. Arguing as before, by the autonomous character of \eqref{CN} we have that $u_2(x,s+t)\le u_2(x,t)$ for any $x\in\mathcal{G}$ and $s\in\R_+$\,. Hence the map  $t\mapsto u_2(\cdot,t)$ is nonincreasing. Similarly, the map $t\mapsto u_1(\cdot,t)$ is nondecreasing.

It remains to prove the statements concerning $\hat{u}_1$ and $\hat{u}_2$. The proof is the same in either case, thus we only address the latter. Let $v\in \mathcal{G}\setminus\pa\mathcal{G}$ be fixed, and let $\Sig_v$ be the star centered at $v$. Let $\phi: \Sig_v\times \R_+\mapsto\R_+$\,,  $\phi(x,t):=\sum_{e\ni v}\phi _e(x)\theta _e(t)$ satisfy the following assumptions:

\noindent $(a)$ $\phi _e\in C^2(I _e)\cap C^1(\overline{I}_e)$, $\phi _e\ge0$, with compact support in $[0,l_e)$  if $e\subseteq\Sig_v^-$, respectively in $(0,l_e]$ if $e\subseteq\Sig_v^+$;

\noindent $(b)$ $\theta _e\in C^1(\R_+)\cap L^\infty(\R_+)$, $\theta _e\ge0\,,$ $\theta_e'\ge0$\,.

 Let $\t\in \R_+$ and $e\ni v$ be fixed. Multiplying \eqref{solge1} (written for $u=u_2$) by $\theta _e\phi _e$ and integrating  in $
I_e\times(0,\t)$ gives
\begin{equation}\label{bar1}
0\,=\, \int_0^\t dt\, \theta _e(t)\int_0^{l_e} dx\left[u_{2et}(x,t)-u_{2exx}(x,t) - f(u_{2e}(x,t))\right]\phi_e(x)\,.
\end{equation}
Plainly, there holds
\begin{eqnarray*}
&&\int_0^\t dt\,\theta _e(t)\!\!\int_0^{l_e}dx\,u_{2et}(x,t)\phi_e(x)\,=\\
&=&
\theta _e(\t)\!\!\int_0^{l_e} u_{2e}(x,\t)\phi _e(x)\,dx \,-\, \theta _e(0)\int_0^{l_e} \overline{q}(x)\phi _e(x)\,dx \,-\!
\int_0^\t \!dt \,\theta _e'(t)\!\!\int_0^{l_e} u_{2e}(x,t)\phi _e(x)\,dx \,, \nonumber
\end{eqnarray*}
and for any $e\subseteq\Sig_v^\mp$
\begin{eqnarray*}
\int_0^\t dt\,\theta _e(t)\!\!\int_0^{l_e} dx\,u_{2exx}(x,t)\phi_e(x)&=& \mp \phi_e(v) \int_0^\t u_{2ex}(v,t)\theta _e(t)\,dt
\pm \phi_e'(v) \int_0^\t u_{2e}(v,t)\theta _e(t)\,dt\,+\\
&+&\int_0^\t dt\,\theta _e(t)\!\!\int_0^{l_e} dx\,u_{2e}(x,t)\phi_e''(x)\, .  \nonumber
\end{eqnarray*}
Then from  \eqref{bar1} for any $e\subseteq\Sig_v^\mp$ we obtain
\begin{eqnarray*}
&&\frac1\t\;\Big\{\,\theta _e(\t)\!\!\int_0^{l_e} u_{2e}(x,\t)\phi _e(x)\,dx \,-\, \theta _e(0)\int_0^{l_e} \overline{q}(x)\phi _e(x)\,dx \,-\!
\int_0^\t \!ds\, \theta _e'(t)\!\!\int_0^{l_e} u_{2e}(x,t)\phi _e(x)\,dx \Big\}=\\
&=& \frac1\t\;\Big\{\mp \phi_e(v) \int_0^\t u_{2ex}(v,t)\theta _e(t)\,dt
\pm \phi_e'(v) \int_0^\t u_{2e}(v,t)\theta _e(t)\,dt\
\,+\\
&+&\int_0^\t \! dt\,\theta _e(t)\int_0^{l_e}\big[ u_{2e}(x,t)\phi _e''(x)+f(u_{2e}(x,t))\phi _e(x)\big]\,dx
 \Big\}\,.    \nonumber
\end{eqnarray*}
Letting $\t\to\infty$ in the above equality, by de $\rm{l'H\hat{o}pital's}$ rule we get plainly
\begin{eqnarray}\label{bar2}
&&
\pm \phi_e(v) \lim_{\t\to\infty}\frac1\t\int_0^\t u_{2ex}(v,t)\theta _e(t)\,dt\,=\\
&=& \theta_e^\infty\,\Big\{\pm\phi_e'(v)\,\hat{u}_2(v) \,+\, \int_0^{l_e}\big[ \hat{u}_{2e}(x)\phi _e''(x)
+f(\hat{u}_{2e}(x))\phi _e(x)\big]\,dx\Big\}  \qquad(e\subseteq\Sig_v^\mp)   \nonumber
\end{eqnarray}
for any $\phi_e$ as above (observe that the limits $\lim_{t\to\infty}\theta_e(t)=:\theta_e^\infty\in\R_+$\,, $\lim_{t\to\infty}\theta_e'(t)=0$ and
$\lim_{t\to\infty}u_{2e}(\cdot,t)=:\hat{u}_{2e}$ exist by the monotone character of $\theta_e$ and of the map $t\mapsto u_{2e}(\cdot,t)$).

From \eqref{bar2} we get
$$
\int_0^{l_e}\big[ \hat{u}_{2e}(x)\phi _e''(x) +f(\hat{u}_{2e}(x))\phi _e(x)\big]\,dx\,=\,0 \quad\text{for any $\phi_e\in C_c^2(I_e)$ and $e\subseteq\Sig_v$\,,}
$$
whence by standard arguments
\begin{equation}\label{stazn}
\hat{u}_{2e}'' +\, f(\hat{u}_{2e})\,=\,0 \quad\text{ in $I_e$} \qquad (e\subseteq\Sig_v)\,.
\end{equation}

Multiplying equality \eqref{stazn} by $\phi _e$ and integrating in $I_e$ gives
\begin{equation*}
\pm \phi_e(v) \hat{u}_{2e}'(v) \,=\,
\pm\phi_e'(v)\,\hat{u}(v) \,+\, \int_0^{l_e}\big[ \hat{u}_{2e}(x)\phi _e''(x) +f(\hat{u}_{2e}(x))\phi _e(x)\big]\,dx \qquad(e\subseteq\Sig_v^\mp) \,,
\end{equation*}
whence by \eqref{bar2} we obtain
\begin{equation}\label{bar3}
 \lim_{\t\to\infty}\frac1\t\int_0^\t u_{2ex}(v,t)\theta _e(t)\,dt\,=\,\theta_e^\infty\, \hat{u}_{2e}'(v)\quad \text{for all $e\subseteq\Sig_v$}\,.
\end{equation}

Since $u_2=\bigoplus_{e\in E}u_{2e}$ is a solution of \eqref{CN}, from \eqref{bar3} with $\theta_e=1$ for all $e\subseteq\Sig_v$ and \eqref{tolex1} we get
\begin{eqnarray*}
0&=&
 \lim_{\t\to\infty}\frac1\t\int_0^\t dt\;\Big\{\sum_{k=1}^{d_v^+} u_{2{e_k}x}(v,t)\,-\,\sum_{l=1}^{d_v^-} u_{2{e_l}x}(v,t) \Big\}\,=\\
 &=& \sum_{k=1}^{d_v^+} \hat{u}_{2{e_k}}'(v)\,-\,\sum_{l=1}^{d_v^-} \hat{u}_{2{e_l}}'(v)\,=\, \sum_{e\ni v}\frac{d \hat{u}_{2e}}{d\n}(v) \quad \text{for all $v\in  \mathcal{G}\setminus\pa\mathcal{G}$}\,,
\nonumber
\end{eqnarray*}
where $\hat{u}_2:=\bigoplus_{e\in E}\hat{u}_{2e}$.

\smallskip

Therefore, $\hat{u}_2$ satisfies the Kirchhoff condition at any $v\in  \mathcal{G}\setminus\pa\mathcal{G}$. It is immediately seen that, if $\pa\mathcal{G}\neq\emptyset$, it also satisfies the boundary conditions at $\pa\mathcal{G}$ both in the Neumann and in the Dirichlet case, thus by \eqref{stazn} it is a stationary solution of \eqref{diffeq}. The maximality of $\hat{u}_2$ in the interval $[\underline{q},\overline{q})$ follows from Theorem \ref{geco} by a standard argument. Hence the result follows.
\hfill$\square$

\smallskip

\rosso{To prove Theorem \ref{comparison}} some preliminary remarks are in order. Let $ \lambda >0$ and $\t\in(0,T]$ be fixed. Let $\z=\bigoplus_{e\in E}  \z_e\in C_0^\infty(\mathcal{G})$, $0\le \z\le 1$\,, and let $n_0\in \mathbb N$ be so large  that $\operatorname{supp}\z\subset \mathcal G_{n_0}$. Then for any $n\in \mathbb N, n\ge n_0+1$ there exists a unique solution $\varphi_n =\bigoplus_{e\in E}  \varphi_{ne}$ to  the backward problem:
\begin{equation*}
\begin{cases}
\varphi_{nt} = -\D \varphi_n +  \lambda \varphi_n &\textrm{in}\,\,\mathcal G_n\times (0, \tau)
\smallskip\\
\varphi_n = 0 & \text{in } (S_n\setminus\pa\mathcal{G})\times (0, \tau)
\smallskip\\
\varphi_n = \z
& \textrm{in}\,\, \mathcal{G}_n\times \{\t\} \,.
\end{cases}
\leqno{(BP_n')}
\end{equation*}
If $\pa\mathcal{G}\neq \emptyset$, since $\varphi_n(\cdot,t)\in D(\D)$, for any $t\in(0,\t)$ there holds:
\begin{equation}\label{duffop4n}
\text{for all $v\in  \pa\mathcal{G}$ }
\quad\left\{
\begin{array}{ll}
\varphi_{nx}(v,t)\,=\,0 \;\;\text{in the Neumann case,}
\\& \\
\varphi_n(v,t)\,=\,0 \;\;\text{in the Dirichlet case}.
\end{array}
\right.
\end{equation}
Moreover, by the maximum principle
\begin{equation}\label{e5f}
0\leq \varphi_n(x,t) \,\le\, e^{\lambda(t-\tau)} \quad \text{ for any } (x,t)\in\mathcal G_n\times (0, \tau)\,.
\end{equation}
Observe that by \eqref{e1f2}, and by \eqref{duffop4n} and the first inequality in \eqref{e5f} if $\pa\mathcal{G}\neq \emptyset$, there holds
\begin{equation}\label{n43}
\sum_{e\in E_n, e\ni v}\frac{d \varphi_{ne}}{d\n}(v,t)\leq 0 \quad \textrm{for any } v\in S_n\times (0, \tau)\,.
\end{equation}

Now we can state a preliminary result, whose proof is given at the end of the section.
\begin{lem}\label{lemma2}
Let the assumptions of Theorem \ref{comparison} be satisfied. Then
\begin{equation}\label{e25f}
\lim_{n\to \infty} \, \sup_{t\in [0, \tau]}\;\sum_{v\in S_n, e\in E_n, e\ni v}\frac{d \varphi_{ne}}{d\n}(v,t) \,=\, 0\quad\text{for any $\tau\in(0,T)$ sufficiently small}\,.
\end{equation}
\end{lem}

Relying on Lemma \ref{lemma2} we can prove the following proposition.
\begin{prop}\label{stima1}
Let $\mathcal G$ be infinite, and let $(H_0)$-$(H_1)$ hold. Let
\begin{equation}\label{e25f-}
\lim_{n\to \infty} \, \sup_{t\in [0, \tau]}\;\sum_{v\in S_n\setminus\pa\mathcal G, e\in E_n, e\ni v}\frac{d \varphi_{ne}}{d\n}(v,t) \,=\, 0\quad\text{for any $\tau\in(0,T)$ sufficiently small}\,,
\end{equation}
and let  $\underline u$, $\overline u$ be a bounded subsolution, respectively a bounded supersolution to \eqref{CN}. Then for any $\lambda>0$ and $t\in [0,T]$ there holds
\begin{equation}\label{e2f}
e^{\lambda t}\int_{\mathcal G}\,[\underline u(x,t)-\overline u(x,t)]_+\, dx\,\le\, \int_0^t \!ds\,e^{\lambda s}\!\!\int_{\mathcal G} \,[f(\underline u)- f(\overline u)+\lambda(\underline u-\overline u)]_+\, dx\,.
\end{equation}
\end{prop}
\begin{proof} Let $\tau\in(0,T)$ be so small that \eqref{e25f} holds, and let $n \ge n_0+1$ be arbitrarily fixed. As in the proof of Theorem \ref{geco}, for all $e\in E_n$ there holds
\begin{equation*}
w_{et} - w_{exx} \,\le\,  f(\underline u_e) - f(\overline u_e)\quad \text{ for any } x\in I_e, t\in (0, \tau]\,,
\end{equation*}
where $w:=\underline{u}-\overline{u}=\bigoplus_{e\in E} w_e$. Multiplying the above inequality by $\varphi_n$ and integrating by parts plainly gives:
\begin{eqnarray}\label{e14f}
&&\int_0^\t \!\!dt\!\int_{\mathcal G_n}[f(\underline u)-f(\overline u)]\,\varphi_n \,dx \,\ge \\
&\ge&\int_{\mathcal G_n} w(x,\tau)\,\z(x)\, dx -\int_{\mathcal G_n} w(x,0)\varphi_n(x,0)\, dx -\lambda\int_0^\t\!\!dt\!\int_{\mathcal G_n}  w\,\varphi_n\, dx\,- \nonumber\\
&-&  \sum_{v\in V_n \setminus S_n}\int_0^\t\Big[ \sum_{k=1}^{d_{v,n}^+}w_{{e_k}x}(j(e_k),t)
\,-\, \sum_{l=1}^{d_{v,n}^-}w_{{e_l}x}(i(e_l),t)\Big]\varphi_n(v,t)\,dt  \,+\nonumber\\
&+&  \sum_{v\in  V_n \setminus S_n}\int_0^\t
\Big[ \sum_{k=1}^{d_{v,n}^+}\varphi_{n{e_k}x}(j(e_k),t) \,-\, \sum_{l=1}^{d_{v,n}^-}\varphi_{n{e_l}x}(i(e_l),t) \Big]w(v,t)\,dt  \,+\nonumber\\
&+&\sum_{v\in S_n}\int_0^{\tau} \Big[w(v, t) \sum_{k=1}^{d_{v,n}^+}\varphi_{n{e_k}x}(j(e_k),t) -\varphi_n(v, t)  \sum_{k=1}^{d_{v,n}^+}w_{{e_k}x}(j(e_k),t)\Big]\, dt\,= \nonumber\\
&=& \int_{\mathcal G_n} w(x,\tau)\z(x) dx -\int_{\mathcal G_n} w(x,0)\varphi_n(x,0) dx -\lambda\int_0^\t\!\!dt\!\int_{\mathcal G_n}  w\,\varphi_n\, dx\,+ \nonumber\\
&+&\sum_{v\in V_n, e\in E_n, e\ni v }\int_0^\t\Big[\frac{d\varphi_{ne}}{d\n}(v,t)w(v,t)-\frac{dw_e}{d\n}(v,t)\varphi_n(v,t)\Big]\,dt\,;\nonumber
\end{eqnarray}
here use of \eqref{e1f} has been made. Since $w(x,0)\varphi_n(x,0)\,\le\,0$  for any $x\in \mathcal G_n$\,, from \eqref{e14f} we get
\begin{equation}\label{e14z}
\int_0^\t \!\!dt\!\int_{\mathcal G_n}[f(\underline u)-f(\overline u)]\,\varphi_n \,dx \,\ge
\int_{\mathcal G_n} w(x,\tau)\z(x) dx  - \lambda\int_0^\t\!\!dt\!\int_{\mathcal G_n}  w\,\varphi_n\, dx  \,+ \,\sum_{k=1}^3 I_{k,n}\,,
\end{equation}
where
$$
I_{1,n}\,:=\!\sum_{v\in V_n\setminus S_n, e\in E_n, e\ni v }\int_0^\t\Big[\frac{d\varphi_{ne}}{d\n}(v,t)w(v,t)-\frac{dw_e}{d\n}(v,t)\varphi_n(v,t)\Big]\,dt\,,
$$
$$
I_{2,n}\,:=\!\sum_{v\in S_n\cap \pa\mathcal G, e\in E_n, e\ni v }\int_0^\t\Big[\frac{d\varphi_{ne}}{d\n}(v,t)w(v,t)-\frac{dw_e}{d\n}(v,t)\varphi_n(v,t)\Big]\,dt\,,
$$
$$
I_{3,n}\,:=\!\sum_{v\in S_n\setminus \pa\mathcal G, e\in E_n, e\ni v }\int_0^\t\Big[\frac{d\varphi_{ne}}{d\n}(v,t)w(v,t)-\frac{dw_e}{d\n}(v,t)\varphi_n(v,t)\Big]\,dt\,.
$$

Now observe that:

\noindent $(a)$ by Definition \ref{sosol2}, inequality \eqref{solprop23} and the  first inequality in \eqref{e5f} there holds \begin{equation}\label{e16f}
\sum_{v\in V_n\setminus S_n, e\in E_n, e\ni v}\frac{dw_e}{d\n}(v,t)\varphi_n(v,t)\,\le\,0\,;
\end{equation}

\noindent $(b)$  for any $t \in(0,\t)$, since $\varphi_n(\cdot,t)\in D(\D)$ there holds
\begin{equation}\label{e15f}
\sum_{v\in V_n\setminus S_n, e\in E_n, e\ni v }  \frac{d\varphi_{ne}}{d\n}(v,t)\,=\,0\,.
\end{equation}
From \eqref{e16f}-\eqref{e15f} we obtain that $I_{1,n}\ge 0$ $(n\in\N)$.

If  $\pa\mathcal G=\emptyset$, we have that $I_{2,n}=0$ (see Remark \ref{mark}). Otherwise, suppose first that the Dirichlet condition \eqref{dico} holds. Then by \eqref{duffop4n} there holds
\begin{equation*}
I_{2,n}\,=\! \sum_{v\in S_n\cap \pa\mathcal G, e\in E_n, e\ni v }\int_0^\t \frac{d\varphi_{ne}}{d\n}(v,t)w(v,t)\,dt\, . \end{equation*}
Moreover, by Definition \ref{sosol2} and inequalities \eqref{solprop3} and \eqref{n43},
\begin{equation*}
\sum_{v\in S_n\cap \pa\mathcal G, e\in E_n, e\ni v }\int_0^\t \frac{d\varphi_{ne}}{d\n}(v,t)w(v,t)\,dt\,= \, \sum_{v\in S_n\cap\pa\mathcal G}\int_0^{\tau} \varphi_{nx}(v,t)\,w(v, t)\,dt\,\ge\,0\,
\end{equation*}
(recall that $d_{v,n}^+=d_v^+=1$ for any $v\in S_n\cap \pa\mathcal G$). It follows that $I_{2,n}\ge 0$ $(n\in\N)$ in this case. Similarly, if $\pa\mathcal G\neq\emptyset$ and the Neumann condition \eqref{neco} holds, by \eqref{duffop4n} there holds
\begin{equation*}
I_{2,n}\,=\, -\!\sum_{v\in S_n\cap \pa\mathcal G, e\in E_n, e\ni v }\int_0^\t \frac{d w_e}{d\n}(v,t)\varphi_n(v,t)\,dt\,,
\end{equation*}
whereas  by Definition \ref{sosol2} and inequalities \eqref{solprop3} and \eqref{e5f}
\begin{equation*}
\sum_{v\in S_n\cap \pa\mathcal G, e\in E_n, e\ni v }\int_0^\t \frac{d w_e}{d\n}(v,t)\varphi_n(v,t)\,dt\,= \, \sum_{v\in S_n\cap\pa\mathcal G}\int_0^{\tau} w_x(v, t)\,\varphi_n(v,t)\,dt\,\le\,0\,.
\end{equation*}
Therefore, there holds $I_{2,n}\ge 0$ $(n\in\N)$ in this case, too.

Finally, since $\varphi_n=0$  in $ (S_n\setminus\pa\mathcal{G})\times (0, \tau)$, there holds
$$
I_{3,n}\,=\!\sum_{v\in S_n\setminus \pa\mathcal G, e\in E_n, e\ni v }\int_0^\t\frac{d\varphi_{ne}}{d\n}(v,t)\,w(v,t)\,dt\,.
$$

In view of the above remarks, from inequality \eqref{e14z} we get
\begin{eqnarray*}
&&\int_{\mathcal G_n}  w(x,\t)\zeta(x) dx - \lambda\int_0^\t\!\!dt\!\int_{\mathcal G_n}  w\,\varphi_n\, dx \,\le
\\
&\le& \int_0^\t\!\!dt\!\int_{\mathcal G_n}[f(\underline u)-f(\overline u)]\,\varphi_n \,dx  \;-\!\! \sum_{v\in S_n\setminus \pa \mathcal G, e\in E_n, e\ni v }\int_0^\t \frac{d\varphi_{ne}}{d\n}(v,t)w(v,t)\,dt\,,
\end{eqnarray*}
whence by the second inequality in \eqref{e5f}
\begin{eqnarray*}
\int_{\mathcal G_n}  w(x,\t)\zeta(x) dx
 &\le& \int_0^\t\!\!dt\!\int_{\mathcal G_n}[f(\underline u)-f(\overline u)+\lambda  w]_+e^{\lambda(t-\tau)} dx  \;-\\
 &-&\sum_{v\in S_n
\setminus \pa \mathcal G, e\in E_n, e\ni v }\int_0^\t \frac{d\varphi_{ne}}{d\n}(v,t)w(v,t)dt
\end{eqnarray*}
(here $[r]_+:=\max\{r,0\}; r\in\R$). Since $\bigcup_{n=1}^{\infty}\mathcal G_{n}=\mathcal G$, letting $n\to\infty$ in the above inequality and using \eqref{e25f-} gives
\begin{equation}\label{e19f}
e^{\l \tau}\!\!\int_{\mathcal G}  w(x,\t)\,\zeta(x)\, dx\,\le\,  \int_0^\t\!\!dt \,e^{\lambda t}\!\int_{\mathcal G}[f(\underline u)-f(\overline u)+\lambda  w]_+\, dx \,.
\end{equation}
Arguing as in the proof of Theorem \ref{geco} shows that inequality \eqref{e19f} still holds if $\zeta(x)$ is replaced by $ [w(x,\t)]_+\,.$ Then inequality  \eqref{e2f} follows directly if $\t=T$, or  otherwise by iterating the above procedure. This completes the proof.
\end{proof}
\smallskip

Theorem \ref{comparison} follows from Lemma \ref{lemma2} and  Proposition \ref{stima1} by a standard argument (e.g., see \cite[Theorem 12]{ACP}). We give the proof for convenience of the reader.

\smallskip

\noindent {\em Proof of Theorem \ref{comparison}.}
Under the present assumptions \eqref{e25f} holds (see Lemma \ref{lemma2}), thus  Proposition \ref{stima1} can be used. In particular, since $\underline{u},\overline{u}\in L^\infty(\O_T)$, inequality \eqref{LM} holds for all $(x,t)\in \O_T$ with $M:=\max\{\|\underline{u}\|_{L^\infty(\O_T)}\,,\|\overline{u}\|_{L^\infty(\O_T)}\}$. Choosing $\l= L_M$ gives
$$
[f(\underline u)- f(\overline u)+\lambda(\underline u-\overline u)]_+\,\le\, 2L_M\,[\underline u-\overline u]_+\,.
$$
From \eqref{e2f} and the above inequality by Gronwall's lemma the result follows.
\hfill$\square$
\smallskip

Let us finally prove Lemma \ref{lemma2}.
\smallskip

\noindent{\em Proof of Lemma \ref{lemma2}.}
 Let $\z\in C_0^\infty(\mathcal{G})$ as in $(BP_n')$\,, and let $n_0\in \mathbb N$ be so large  that: $(a)$ $\operatorname{supp}\z\subset \mathcal G_{n_0}$\,; $(b)$ $\pa\mathcal G\subset \mathcal G_{n_0}$\,, if assumption $(ii)$ holds. Let $n\ge n_0+1$ and $\t\in(0,T]$ be fixed. Recall that $\mathcal G_n\setminus \mathcal G_{n_0}=\{x\in\mathcal G \,|\,R_{n_0}<r(x)\le R_n \}$, and similarly $\mathcal G_n\setminus \mathcal G_{n-1}=\{x\in\mathcal G \,|\,R_{n-1}<r(x)\le R_n \}$. Consider the function
\begin{equation}\label{defieta}
\eta: \left(\mathcal G_n\setminus \mathcal G_{n_0}\right)\times(0, \tau] \mapsto\R_+\,, \quad
 \eta(x, t):=\sigma \exp\left\{\frac K{t-\tau-t_0}[r(x)+r_0]^\b\right\} \quad  (x\in  \mathcal G_n\setminus \mathcal G_{n_0}, t\in (0, \tau])
\end{equation}
with $\b\in[0,2]$ as in $(H_2)$-$(ii)$, and $\sigma:=\exp{\left\{2\theta(R_{n_0}+r_0)^\b\right\}}$, $t_0:=\tfrac{K}{2\theta}$ with
\begin{equation}\label{sce}
\begin{cases}
K>0\,, \quad r_0:=\max\left\{\left(K\b^2\right)^{\frac{1}{2-\b}}-R_{n_0}, c_0 \right\} &\textrm{if}\,\,\b\in[0,2)\,,
\smallskip\\
K\in\left(0,\frac14\right)\,, r_0:=c_0 & \text{if}\; \b=2\,.
\end{cases}
\end{equation}

Set also
\begin{equation}\label{defih}
h: \mathcal G_n\setminus \mathcal G_{n-1} \mapsto\R\,, \quad
h(x):=\frac{\eta(R_{n-1},0)}{R_n-R_{n-1}}\,[R_n-r(x)] \quad  (x\in  \mathcal G_n\setminus \mathcal G_{n-1})
\end{equation}
with $\eta$ given by \eqref{defieta}.  By the very definition of $h$ we have that
\begin{equation}\label{n35}
h \,=\, 0  \quad \textrm{in}\  S_n\times(0,\t)\,.
\end{equation}
Furthermore, since $h'=-\frac{\eta(R_{n-1},0)}{R_n-R_{n-1}}\,<0$, by \eqref{e1f1} and assumption $(H_2)$-$(i)$ we get
\begin{equation}\label{n40a}
\sum_{e\in E_n, e\ni v}\frac{d h_e}{d\n}(v) =-\,(d^+_v-d^-_v)\,\frac{\eta(R_{n-1},0)}{R_n-R_{n-1}}\,\ge\,0\quad \text{ for all }   v\in V_n\setminus S_n\,,
\end{equation}
whereas by \eqref{e1f2} there holds
\begin{equation}\label{e2lu}
\sum_{e\in E_n, e\ni v}\frac{d h_e}{d\n}(v)\,=\,-\,d^+_v\,\frac{\eta(R_{n-1},0)}{R_n-R_{n-1}}\,<\,0\quad \textrm{for all } v\in S_n\,.
\end{equation}
Similarly, since $\eta(\cdot,t)$ $(t\in(0,\t)$ only depends on the distance $r(x)$, and
\begin{equation}\label{n23}
\eta_x(x, t)=\frac{ K \b (r(x)+r_0)^{\b-1}}{t-\tau-t_ 0}\,\eta(x, t) < 0\,,
\end{equation}
by \eqref{e1f1} and $(H_2)$-$(i)$ we have that
\begin{equation}\label{e22f}
\sum_{e\in E_n, e\ni v}\frac{d \eta_e}{d\n}(v) =-\,(d^+_v-d^-_v)\,
\frac{ K \b (r(v)+r_0)^{\b-1}}{t-\tau-t_0}\,\eta(v, t) \,\ge\,0\quad \text{ for all }   (v,t)\in (V_n\setminus S_n)\times(0,\t)\,.
\end{equation}

We shall prove the following
\smallskip

\noindent {\em  Claim:} The function $h$ is a supersolution of the backward problem
\begin{equation}\label{n41}
\begin{cases}
z_t  +  \Delta z = 0 & \textrm{in } (\mathcal G_n\setminus \mathcal G_{n-1}) \times (0, \tau) \, , \\
z = 0 & \textrm{on } S_n\times (0, \tau) \, , \\
z = \varphi_n & \textrm{on }  S_{n-1} \times (0, \tau) \, , \\
z = 0  & \textrm{in }(\mathcal G_n\setminus \mathcal G_{n-1})\times \{\tau\}\,.
\end{cases}
\end{equation}

From the above Claim the result follows. In fact, since $\varphi_n$ is a solution of $(BP_n')$, it is a subsolution of \eqref{n41} (observe that ${\rm supp}\,\z\subseteq \mathcal G_{n_0}\subseteq \mathcal G_{n-1}$ since $n\ge n_0+1$, thus  $\varphi_n=0$ in $(\mathcal G_n\setminus \mathcal G_{n-1})\times \{\tau\}$).  Then by comparison results there holds
\begin{equation}\label{n42}
 \varphi_n\,\le\, h  \quad \textrm{in } (\mathcal G_n\setminus \mathcal G_{n-1})\times (0, \tau)\,.
\end{equation}

By \eqref{n35} and \eqref{n42}, for any $(v,t)\in S_n\times (0, \tau)$ there holds
\begin{equation*}
\sum_{e\in E_n, e\ni v}\frac{d [h_e-\varphi_{ne}(\cdot,t)]}{d\n}(v)\,\le\,0\,, \end{equation*}
whence by \eqref{n43}
\begin{equation}\label{n45}
\left| \,\sum_{e\in E_n, e\ni v}\frac{d \varphi_{ne}}{d\n}(v,t) \right| \;\le\; \left|\, \sum_{e\in E_n, e\ni v}\frac{d h_e}{d\n}(v) \right| \,.
\end{equation}
In view of \eqref{e2lu} and \eqref{n45}, for any $n\ge n_0+1$ and $(v,t)\in S_n\times (0, \tau)$ there holds
\begin{equation}\label{n46}
\left|\, \sum_{e\in E_n, e\ni v}\frac{d \varphi_{ne}}{d\n}(v,t)\, \right| \,\le  \,d^+_v\,\frac{\eta(R_{n-1},0)}{R_n-R_{n-1}}\,=  \,
d^+_v\, \s \, \exp\left\{-\frac{K}{\tau+t_0}(  R_{n -1}+ r_0)^{\b}\right\}\,.
\end{equation}
From \eqref{e30f}, \eqref{n46} and assumption $(H_2)$-$(ii)$ we obtain that
\begin{eqnarray*}
&&\left|\, \sup_{t\in (0, \tau]}\sum_{v\in S_n, e\in E_n, e\ni v}\frac{d \varphi_{ne}}{d\n}(v,t) \right|\;\le\;
\left(\sum_{v\in S_n}d_v^+\right) \s \, \exp\left\{-\frac{K}{\tau+t_0}(R_{n -1}+ r_0)^{\b}\right\}
\,\le \\
&\le& C \sigma \exp\left\{\left[\theta  -\frac{K}{\tau+t_0}\left(\frac{  R_{n -1}+ r_0}{  R_n}\right)^{\b} \right]  R_n^{\beta} \right\}\,\le\, C \sigma \exp\left\{\theta(\t-t_0) R_n^{\beta}\right\}\,,
\end{eqnarray*}
since $\frac{R_{n -1}+ r_0}{  R_n}\ge\frac{  R_{n -1}+ r_0}{  R_{n -1} +c_0}\ge1$\,. Then for any $\tau\in(0,t_0)$ letting $n\to\infty$ gives \eqref{e25f}.
\medskip

It remains to prove the Claim. To this purpose, let us first prove that $\eta$ is a supersolution of the backward problem
 \begin{equation}\label{n9b}
\begin{cases}
z_t  + \, \Delta z = 0 & \textrm{in }\ (\mathcal G_n \setminus \mathcal G_{n_0}) \times (0, \tau) \, , \\
z = 0 & \textrm{on } S_n\times (0, \tau) \, , \\
z = \varphi_n & \textrm{on } S_{n_0} \times (0, \tau) \, ,\\
z = 0 & \textrm{in } ( \mathcal G_n \setminus \mathcal G_{n_0}) \times \{\tau\} \,.
\end{cases}
\end{equation}
 In fact, for every $(x,t)\in (\mathcal G_n \setminus \mathcal G_{n_0}) \times (0, \tau)$ there holds
 \begin{equation}\label{n22}
\eta_t(x, t)=-\frac{K [r(x)+r_0]^{\b}}{(t-\tau-t_0)^2} \, \eta(x, t) < 0\,,
\end{equation}
\begin{equation*}
\eta_{xx}(x, t)=\left\{\frac{ K^2 \b^2 [r(x)+r_0]^{2\b-2}}{(t-\tau-t_0)^2} + \frac{ K \b(\b-1)[r(x)+r_0]^{\b-2}}{t-\tau-t_0}\right\}\eta(x, t)\,\le\, \frac{ K^2 \b^2 [r(x)+r_0]^{2\b-2}}{(t-\tau-t_0)^2} \,\eta(x, t)\,,
\end{equation*}
whence by \eqref{sce}
\begin{equation}\label{n33}
\eta_t(x, t) + \eta_{xx}(x, t) \,\le\, \eta(x,t)\,\frac{K[r(x)+r_0]^{2\b-2}}{(t-\tau-t_0)^2} \Big\{-[r(x)+r_0]^{2-\b} + K\b^2
\Big\}\,\le\,0\,.
\end{equation}
In addition, by \eqref{e5f} and the definition of $\sigma$ there holds
\begin{equation}\label{n34}
 \varphi_n\,\le\, \eta \quad  \textrm{on } S_{n_0} \times (0, \tau)\,.
 \end{equation}

In view of \eqref{e22f}, \eqref{n33} and \eqref{n34}, $\eta$ is in fact a supersolution of problem \eqref{n9b}. On the other hand, it is easily seen that $\varphi_n$ is a subsolution of the same problem. Then by comparison results there holds
\begin{equation}\label{n38}
\varphi_n\,\le\,\eta  \quad \textrm{in } (\mathcal G_n\setminus \mathcal G_{n_0})\times (0, \tau)\,.
\end{equation}

 By \eqref{n22} and \eqref{n38} there holds in particular
\begin{equation*}\label{n36}
h(x)\,=\,\eta(  R_{n -1},0\,\ge\, \eta(  R_{n -1}, t) \geq \varphi_n\quad \textrm{for all } (x,t)\in S_{n-1}\times (0, \tau]\,.
\end{equation*}
From \eqref{defih}-\eqref{n40a} and the above inequality the Claim follows. This completes the proof.
\hfill$\square$


\section{An example: metric trees}\label{semetre}

When applying Theorems \ref{geco}-\ref{comparison}, it is sometimes possible to take advantage of specific geometrical properties of the graph $\mathcal{G}$ we are dealing with. Let us discuss this issue in the case of regular metric trees (see Subsection \ref{metre}).

Let $\mathcal{T}$ be a regular tree of infinite height. The boundary $\pa\mathcal{T}$ is reduced to one element, the root $O$. For any vertex $v$ of the $n$-th generation, thus belonging to $\mathcal{T}\setminus O$, the inbound star $\Sig_v^+$  centered at $v$ consists of one edge $e$ (i.e., $d_v^+=1$), whereas the outbound star $\Sig_v^-$  consists of $b_n$ edges $e_l$ (i.e., $d_v^-=b_n$). Accordingly, the Kirchhoff condition satisfied at $v$ by a solution of problem \eqref{CN} on $\mathcal{T}$ reads (see \eqref{tolex1}):
\begin{equation}\label{tolex111}
u_{ex}(v,t)\,=\,\sum_{l=1}^{{b_n}} u_{{e_l}x}(v,t)  \,.
\end{equation}

It is often convenient to construct {\em symmetric} sub- and supersolutions (in particular, solutions) of problem \eqref{CN} - namely, sub- and supersolutions which only depend on the distance from the root $O$ (see Definition \ref{symm}). Symmetric sub- and supersolutions are of the form $\underline{u}(x,t)=\underline{z}(\r(x),t)$, respectively $\overline{u}(x,t)=\overline{z}(\r(x),t)$ $(t\in[0,T])$, where $\underline{z},\,\overline{z}:[0,\R_+)\times[0,T]\mapsto \R$ and $\r=\r(x):=d(x,O)$ $(x\in\mathcal{T})$.

In view of Definition \ref{evoss} and of the structure of a regular tree, it is immediately seen that $\overline{u}$ is a symmetric supersolution of problem \eqref{CN}  in $[0,T]$ if and only if $\overline{z}$ belongs to $C([0,T];L^2(\R_+;\b))\cap\, C^1((0,T];L^2(\R_+;\b))\cap C((0,T];H^1(\R_+;\b))$ with the following properties:

\noindent $a)$  for any $t\in[0,T]$ there holds $\overline{z}(\cdot,t)=\sum_{n=1}^\infty \overline{z}_n(\cdot,t)\chi_{_{I_n}}$, where $I_n:=(\r_{n-1},\r_n)$; here $\r_n$ $(n\in\N)$ is the distance from the root $O$ of any vertex of the $n$-th generation and $\r_0:=0$ (see Definition \ref{detree}-$(iii)$);

\noindent $b)$ there holds
\begin{subequations}\label{sotto}
\begin{equation}\label{sotto20}
 \overline{z}_n(\cdot,t)\in H^2(I_n)\;\textrm{ for all $n\in\N$}\,,\quad \sum_{n\in\N}\int_{I_n}
 \left|\overline{z}_{n\r\r} (\r,t) \right|^2\b(\r)  d\r\,<\,\infty\,,
\end{equation}
\begin{equation}
\overline{z}_t(\cdot,t)\,\ge\,  \overline{z}_{\r\r}(\cdot,t) + f(\overline{z}(\cdot,t))  \quad\text{a.e. in $\R_+$ for any $t\in [0,T]$\,,}
\end{equation}
\begin{equation}\label{sotto2}
\overline{z}(\r(x),0)\,\ge\, u_0 \;\;\text{for a.e. $x\in\R_+$}\,,
\end{equation}
\begin{equation}\label{sotto3}
\overline{z}_{\r}(\r_n^-,t)\,\ge\, b_n\;\,\overline{z}_{\r}(\r_n^+,t)\quad\text{for all $n\in\N$}\,,
\end{equation}
\begin{equation}
\left\{
\begin{array}{ll}
\overline{z}_{\r}(0^+,t) \,\le\,0 \;\;\text{in the Neumann case,}
\\& \\
\overline{z}(0,t)\,\ge\,0 \;\;\text{in the Dirichlet case}.
\end{array}
\right.
\end{equation}
\end{subequations}

Similar considerations hold for symmetric subsolutions of problem \eqref{CN}, as well as for symmetric stationary sub- and supersolutions of equation \eqref{diffeq} - namely, for stationary sub- and supersolutions of equation \eqref{diffeq} of the form  $\underline{q}(x)= \underline{Q}(\r(x))$, $\overline{q}(x)= \overline{Q}(\r(x))$ $(x\in\mathcal{T})$. For instance, in view of Definition \ref{soso}, $\overline{q}$ is a symmetric stationary supersolution of equation \eqref{diffeq} if and only if $\overline{Q}$ belongs to $H^1(\R_+;\b)$,
$\overline{Q}=\sum_{n=1}^\infty \overline{Q}_n\chi_{_{I_n}}$, and the following holds:
\begin{subequations}\label{sponge}
\begin{equation}\label{sponge1}
\overline{Q}_n\in H^2(I_n) \,\textrm{ for all $n\in \N$}, \;
 \quad  \sum_{n\in \N}\int_{l_n}| \overline{Q}_n''(\r)|^2\b(\r)\,d\r<\infty \,,
\end{equation}
\begin{equation}\label{sponge2}
 \overline{Q}''+f(\overline{Q})\,\le\,0  \quad\text{a.e. in $\R_+$\,,}
\end{equation}
\begin{equation}\label{sponge3}
\overline{Q}'(\r_n^-)\,\ge\, b_n\;\,\overline{Q}'(\r_n^+)\quad\text{for all $n\in\N$}\,,
\end{equation}
\begin{equation}\label{sponge4}
\quad\left\{
\begin{array}{ll}
 \overline{Q}'(0^+)\le0 \;\; \text{ in the Neumann case,}
\\& \\
 \overline{Q}(0)\,\ge\,0 \;\;\text{in the Dirichlet case}.
\end{array}
\right.
\end{equation}
\end{subequations}
\smallskip

In view of the above remarks, the following comparison result is a typical application of Theorem \ref{decre} to regular trees
(a similar result for bounded solutions follows from Theorem \ref{comparison}).
\begin{theorem}\label{confronto}
Let $\mathcal{T}$ be a regular tree. Let $(H_1)$ hold, and let $f(0)=0$. Let $\overline{Q}\in H^1(\R_+;\b)$,
$\overline{Q}=\sum_{n=1}^\infty \overline{Q}_n\chi_{_{I_n}}$, $\overline{Q}\ge0$ satisfy \eqref{sponge}. Let $\hat u$ be the solution of problem \eqref{CN}  with initial data $\hat u_0=\overline{Q}\circ\r$, and let $0\le u_0\le \hat u_0$
in $\mathcal{T}$. Then the solution $u$ of problem \eqref{CN} is global, and for all $t\in\R_+$ there holds
$0\le u(\cdot,t)\le  \hat u(\cdot,t)$ in $\mathcal{T}$.
\end{theorem}
\begin{proof}
By the above remarks $\hat u_0$ is a symmetric stationary supersolution of \eqref{diffeq}, whereas $0$ is a (symmetric) stationary subsolution. Then the conclusion follows from Theorem \ref{decre}.
\end{proof}


\section*{Appendix}\label{ar}
\renewcommand{\thesection}{\Alph{section}}\setcounter{section}{1}\setcounter{subsection}{0}\setcounter{equation}{0}

While referring the reader to \cite{BK, Mu} for general definitions and results concerning metric graphs, we recall some general notions for convenience.


\subsection{Metric graphs}\label{megra}

Like combinatorial graphs  (e.g., see \cite{Mu}), a metric graph consists of a countable set of vertices and of a countable set of edges. However, in contrast to combinatorial graphs, the edges are regarded as intervals glued together at the vertices:

\begin{definition}\label{meg} A {\em metric graph}  is a quadruple $(E, V, i, j)$, where:

\noindent $(a)$ $E$  is a countable family of open intervals $I_e\equiv(0,l_e)$, called {\em edges} of {\em length} $l_e\in(0,\infty]$;

\noindent $(b)$ $V$ is a countable set of points, called {\em vertices};

\noindent $(c)$ the maps $ i:E\mapsto V$ and $j:\{e\in E\,|\, l_e<\infty\} \mapsto V$ define the {\em initial point} and the {\em final point} of an edge. Both the initial  and the final point are {\em endpoints} of an edge.

\noindent A metric graph is {\em finite} if both $E$ and $V$ are finite.  A finite metric graph with $l_e<\infty$ for all $e\in E$ is called {\em compact}.
\end{definition}
We always assume for simplicity that $i(e)\neq j(e)$ (namely, no {\em loops} exist), and $l_e<\infty$ for all $e\in E$.
The following notations will be used:
$$
\mathcal{G}_e :=\{e\}\cup  I_e\,, \quad \mathcal{G}:= (\bigcup_{e\in E}\mathcal{G}_e)\cup V\,, \quad
 \overline{\mathcal{G}}_e:=\mathcal{G}_e \cup\{i(e),j(e)\}\,.
$$
By abuse of language, we currently speak of the metric graph $\mathcal{G}$.

Observe that $\mathcal{G}_e$ basically coincides with the interval $I_e$, the component $\{e\}$ in its definition being only added to make the $\mathcal{G}_e$'s mutually disjoint. Therefore, no confusion arises, if the same notation $x,y,\dots$ is used both for points of the edge $e\in E$ and for points of the interval $I_e\subseteq\R_+$\,.
For any $e\in E$ the map $\pi_e : \mathcal{G}_e\mapsto I_e$, $\pi_e (\{e\},x)\equiv\pi_e (x):=x$ establishes a one-to-one correspondence between points of $e$ and points of $I_e$\,, which can be extended to a map from $\overline{\mathcal{G}}_e$ to $\bar{I}_e=[0,l_e]$ such that $\pi_e (i(e)) = 0$ and $\pi_e (j(e)) = l_e$\,. Namely, the coordinate $x\in I_e$ assigned to each point $\pi_e^{-1} (x)\in e$  increases from $x=0$ (which corresponds to $\pi_e^{-1} (0)=i(e)$) to $x=l_e$  (which corresponds to $\pi_e^{-1} (l_e)=j(e)$).

\smallskip

We shall write $e\ni v$ or $v\in e$, if either $i(e)=v$ or $j(e)=v$ $(e\in E, v\in V)$. The notation $e\equiv(u,v)$ means that $i(e)=u$ and $j(e)=v$; therefore, $(u,v)\neq(v,u)$ $(e\in E, \, u,v\in V)$.

\begin{definition}\label{deg}
 Let $\mathcal{G}$ be a metric graph.

  \noindent $(i)$  The {\em degree} $d_{v}\in\N$ of a vertex $v$ is the number of edges $e\ni v$. The {\em inbound degree} $d_{v}^+$ (respectively, {\em outbound degree} $d_{v}^-$) is the number of edges with $j(e)=v$ (respectively, with $i(e)=v$).  $\mathcal{G}$ is {\em locally finite} if $d_v<\infty$ for any $v\in V$.

\noindent $(ii)$  The set $\pa \mathcal{G}:=\{v\in V\;|\; d_v=1\}$ is called the {\em boundary} of the graph. The vertices $v\in V\setminus \pa\mathcal{G}$ are called {\em interior vertices}.

  \noindent $(iii)$ The {\em star centered at} $v$,
\begin{equation}\label{sv2}
 \Sig_v \,:=\,\{x\in\mathcal{G}\;|\; d(x,v)< l_j, \, j=1,\dots,d_v\}\qquad (v\in V)\,.
\end{equation}
We also define the {\em inbound star} $\Sig_v^+$ and the {\em outbound star $\Sig_v^-$ centered at} $v$ as the union of the edges ending at $v$, respectively emanating from $v$:
\begin{equation}\label{sv1}
\Sig_v^\pm \,:=\, \{x\in\mathcal{G}\;|\; d(x,v)< l_j, \, j=1,\dots,d_v^\pm\}\,,
\end{equation}
\end{definition}

Clearly, there holds $d_{v}=d_{v}^++d_{v}^-$ and $\Sig_v =\Sig_v^+(v) \cup \Sig_v^-(v)$  $(v\in V)$.

\begin{definition}\label{sig}
$(i)$  Let $\mathcal{G}$ be a metric graph, and let $u,v\in V$. A {\em path} connecting $u$ and $v$ is a set $\{x_1,\dots,x_n\}\subseteq \mathcal{G}$ $(n\in\N)$ such that  $x_1=u$, $x_n=v$ and for all $k=1,\dots, n -1$ there exists an edge $e_k$ such that $x_k, x_{k+1} \in \overline{\mathcal{G}}_{e_k}$\,. A path connecting $u$ and $v$ is {\em closed} if $u\equiv v$. A closed path is called a {\em cycle}, if it does not pass through the same vertex more than once.

\noindent $(ii)$ A metric graph $\mathcal{G}$ is {\em connected} if for any $u,v\in V$ there exists a path connecting $u$ and $v$. A connected graph without cycles is called a {\em tree}.
\end{definition}

 It is easily seen that a connected metric graph $\mathcal{G}$ can be endowed with a structure of metric measure space. In fact, every two points $x,y\in\mathcal{G}$ can be regarded as vertices of a path $P$ connecting them  (possibly adding them to $V$). The length of the path is the sum of its $n$ edges $e_k$, i.e., $l(P):= \sum_{k=1}^n l_{e_k}\,,$ and the distance $d=d(x,y)$ between $x$ and $y$ is
$$
d(x,y)\,:=\,\inf\,\{l(P)\,|\, \text{$P$ connects $x$ and $y$}\}\,.
$$
Therefore $\mathcal{G}$ is a metric space, thus a topological space endowed with the metric topology.

\smallskip

Let $\mathcal{B}=\mathcal{B}(\mathcal{G})$ be the  Borel $\s$-algebra on $\mathcal{G}$. A Radon measure $\m:\mathcal{B}\mapsto [0,\infty]$ is induced on $\mathcal{G}$ by the Lebesgue measure $\l$ on each interval $I_e$, namely
\begin{equation}\label{demu}
\m(G)\,:=\, \sum_{e\in E} \l(I_e\cap G) \quad\text{for any $G\in\mathcal{B}$}\,.
\end{equation}

Let $\mathcal{G}$ be a metric graph. Every function $f:\mathcal{G}\mapsto\R$ canonically induces a countable family $\{f_e\}$, $f_e:I_e\mapsto\R$ $(e\in E)$. This is expressed as $f=\bigoplus_{e\in E}f_e$\,. We set $f^{(h)}:=\bigoplus_{e\in E}f_e^{(h)}$ $(h\in\N)$ if the derivative $f_e^{(h)}\equiv\frac{d^h f_e}{dx^h}$ exists in $I_e$ for all $e\in E$ (we shall also denote $f^{(0)}\equiv f$, $f^{(1)}\equiv f'$ and $f^{(2)}\equiv f''$). Let $C(\mathcal{G})$ denote the space of continuous functions on the metric space $(\mathcal{G},d)$. We set
$$
C^k(\mathcal{G})\,:=\,\{f\in C(\mathcal{G})\;|\; f_e\in C^k(I_e) \;\forall e\in E, \;\; f^{(h)}\in C(\mathcal{G})\;\forall \;h=1,\dots,k\} \qquad(k\in\N)\,,
$$
and $C^0(\mathcal{G})\equiv C(\mathcal{G})$\,. We also denote by $C^k_c(\mathcal{G})$ the subspace of functions in $C^k(\mathcal{G})$ with compact support, by $C_0(\mathcal{G})$ the closure of $C_c(\mathcal{G})\equiv C_c^0(\mathcal{G})$ with respect to the norm $\|\cdot\|_\infty$, and we set $C^\infty_c(\mathcal{G}):=\bigcap_{k\in\N}C^k_c(\mathcal{G})$\,.

In view of \eqref{demu}, for any measurable $f:\mathcal{G}\mapsto\R$ we set
$$
\int_\mathcal{G} f\,d\m\,:=\, \sum_{e\in E} \int_0^{l_e} f_e\,dx\,, \qquad \int_G f\,d\m\,:=\, \int_\mathcal{G} f\chi_G\,d\m
\quad\text{for any $G\in\mathcal{B}(\mathcal{G})$\,,}
$$
where $\chi_G$ denotes the characteristic function of the set $G$ and the usual notation $dx\equiv d\l$ is used. Accordingly, for any $p\in[1,\infty]$ the Lebesgue spaces $L^p(\mathcal{G})\equiv L^p(\mathcal{G},\m)$ are immediately defined:
$$
L^p(\mathcal{G})\,:=\,\bigoplus_{e\in E}\,L^p((I_e),\l)
$$
with norm
$$
\|f\|_p\,:=\,\sum_{e\in E}\|f_e\|_p\, =\, \sum_{e\in E} \left(\int_0^{l_e} |f_e|^p\,dx\right)^{\frac1p} \;\; \text{if $p\in (1,\infty)$}\,, \quad \|f\|_\infty\,:=\,{\rm ess\, sup}_{x\in\mathcal{G}}\, |f(x)|\,.
$$

For any $p\in[1,\infty]$ and $m\in\N$ the Sobolev spaces are
$$
W^{m,p}(\mathcal{G},\m)\,:=\,\left\{f\in L^p(\mathcal{G},\m)\;|\; f^{(h)}\in L^p(\mathcal{G},\m)\cap C(\mathcal{G})\;\forall \;h=0,\dots,m-1, \; f^{(m)}\in L^p(\mathcal{G},\m)  \right\}
$$
with norm $\|f\|_{m,p}\,:=\,\sum_{h=0}^m\|f^{(h)}\|_p$\,. As usual, we set $W^{m,p}(\mathcal{G})\equiv W^{m,p}(\mathcal{G},\m)$, $W^{0,p}(\mathcal{G})\equiv L^p(\mathcal{G})$,  $H^m(\mathcal{G}):= W^{m,2}(\mathcal{G})$, and we denote by $W^{m,p}_0(\mathcal{G})$ (in particular, by
$H^m_0(\mathcal{G})$) the closure of $C_c^\infty(\mathcal{G})$ with respect to the norm $\|\cdot\|_{m,p}$ (respectively $\|\cdot\|_{m,2}$). If $\mathcal{G}$ is compact, the embedding of $W^{1,p}(\mathcal{G})$ in $C_0(\mathcal{G})$ is continuous for $p\in[1,\infty)$ and compact for $p\in(1,\infty)$ (see \cite[Corollary 2.3]{H}).

For any $V_D\subseteq V$ we define
\begin{equation}\label{dewey}
W^{m,p}(\mathcal{G},V_D)\,:=\,\left\{f\in W^{m,p}(\mathcal{G})\;|\; f^{(h)}(v)=0 \;\forall\,v\in V_D, \;h=0,\dots,m-1 \right\},
\end{equation}
and $H^m(\mathcal{G},V_D)\equiv W^{m,2}(\mathcal{G},V_D)$. Clearly, there holds $W^{m,p}(\mathcal{G},\emptyset)=W^{m,p}(\mathcal{G})$,  and $W^{m,p}(\mathcal{G},\pa \mathcal{G})=W^{m,p}_0(\mathcal{G})$.   In particular, there holds
$$
H^1(\mathcal{G},V_D)\,=\,\left\{f\in H^1(\mathcal{G})\;|\; f(v)=0 \;\forall\,v\in V_D \right\},
$$
thus
$$
H^1(\mathcal{G},\pa\mathcal{G})\,=\,\left\{f\in H^1(\mathcal{G})\;|\;  f(v)=0 \;\forall\,v\in \pa\mathcal{G} \right\}\,=\,H_0^1(\mathcal{G})\,.
$$


\subsection{Laplacian on metric graphs}\label{lamegra}

Using the above functional framework a metric graph $\mathcal{G}$ can be endowed with a Dirichlet form, which makes it a  metric Dirichlet graph. The most obvious choice is the energy form
\begin{equation} \label{nufo}
\mathfrak{a}(f,g)\,:=\, \int_\mathcal{G} f'g'\,d\m\,,
\qquad\text{$f,g\in\mathfrak{D}(\mathfrak{a})\,:=\,H^1(\mathcal{G})$}
\end{equation}
(this choice is feasible, since the Hilbert space $H^1(\mathcal{G})$ is densely defined and continuously embedded in $L^2(\mathcal{G})$). The Laplacian $\D$ associated with $\mathfrak{a}$ is the operator with domain
\begin{subequations} \label{lave}
\begin{equation}
D(\D)\,:=\, \left\{f \in H^1(\mathcal{G})\;|\;
f_e\in H^2(I_e) \,\forall e\in E, \; \sum_{e\in E}\int_0^{l_e}|f_e''|^2\,dx<\infty , \; \sum_{e\ni v}\frac{df_e}{d\n}(v)=0\,\forall\;v\in V\right\}
\end{equation}
(here $\frac{df_e}{d\n}(v)$ denotes the outer derivative of $f_e$ at the vertex $v$), which acts on $D(\D)$ in the natural way:
\begin{equation}
(\D f)_e\,:=\, f_e'' \quad \textrm{for any}\; f\in D(\D) \qquad(e\in E)\,.
\end{equation}
\end{subequations}

Since $H^1(\mathcal{G})\subseteq C(\mathcal{G})$, every $f\in D(\D)$ is continuous at any vertex $v\in \mathcal{G}\setminus  \pa\mathcal{G}$. The condition $\sum_{e\ni v}\frac{df_e}{d\n}(v)=0$ is the Kirchhoff transmission condition if $v\in \mathcal{G}\setminus \pa\mathcal{G}$, and the Neumann homogeneous boundary condition if $v\in \pa\mathcal{G}$\,. Hence the Laplacian defined in \eqref{lave} is called {\em Neumann Laplacian} and denoted by $\Delta_N$, if clarity so requires.

If $V_D\subseteq V$, $V_D\neq\emptyset$, a different notion of Laplacian is obtained considering the Dirichlet form
\begin{equation} \label{difo}
\mathfrak{a}_{_{V_D}}(f,g)\,:=\,\int_\mathcal{G} f'g'\,d\m\,,
\qquad\text{$f,g\in\mathfrak{D}(\mathfrak{a}_{_{V_D}})\,:=\,H^1(\mathcal{G},V_D)$.}
\end{equation}
The Laplacian $\D_{_{V_D}}$ associated with $\mathfrak{a}_{_{V_D}}$ is the operator defined as follows:
\begin{subequations} \label{lavedi}
\begin{equation}
(\D_{_{V_D}} f)_e\,:=\, f_e'' \quad \textrm{for any}\; f\in D(\D_{_{V_D}}) \quad(e\in E)\,,
\end{equation}
where
\begin{equation}
D(\D_{_{V_D}}):= \left\{f \in H^1(\mathcal{G},V_D)\,|\, f_e\in H^2(I_e) \,\forall e\in E,\,\sum_{e\in E}\int_0^{l_e}|f_e''|^2\,dx<\infty , \, \sum_{e\ni v}\frac{df_e}{d\n}(v)=0\,\forall\;v\in V\setminus V_D\right\}.
\end{equation}
\end{subequations}

More explicitly, the vertex conditions satisfied by all $f\in D(\D_{_{V_D}})\subseteq C(\mathcal{G})$ are:
$$
\left\{
\begin{array}{ll}
f(v)=0 &\quad\forall v\in V_D\,, \medskip \\
\sum_{e\ni v}\frac{df_e}{d\n}(v)=0 &\quad\forall v\in V\setminus V_D\ \, .
\end{array}
\right.
$$
Therefore, it is natural to call {\em Dirichlet Laplacian} the Laplacian $\D_D\equiv \D_{ \pa\mathcal{G}}$ (observe that the Neumann Laplacian satisfies $\D_N\equiv \D_\emptyset$). In particular, for a metric tree the vertex conditions satisfied by any $f\in D(\D_D)\subseteq C(\mathcal{G})$ explicitly read:
$$
\left\{
\begin{array}{ll}
f(O)=0\,, & \medskip \\
\sum_{e\ni v}\frac{df_e}{d\n}(v)=0 &\quad\forall v\in \mathcal{G}\setminus O \, .
\end{array}
\right.
$$

\smallskip

Any realization of the Laplacian on a metric graph $\mathcal{G}$ ia a self-adjoint  nonpositive operator, thus the spectrum $\s(-\D)$ is contained in $\overline{\R}_+$.
We shall denote by $\{e^{\D t}\}_{t\ge0}$ the analytic semigroup generated by $\D$ in $L^2(\mathcal{G})$, which exists by classical results (e.g., see \cite{Y}). As in \cite[Theorems 1.4.1-1.4.2]{D}, the following holds:
\begin{prop}\label{davies}
The semigroup $\{e^{\D t}\}_{t\ge0}$ may be extended to a positivity preserving contraction semigroup $\{T_p(t)\}_{t\ge0}$ on $L^p(\mathcal{G})\equiv L^p(\mathcal{G},\m)$ for all $p \in[1,\infty]$. Moreover,  there holds
\begin{equation}\label{consi}
T_p(t)u= T_q(t)u \quad\text{for all $p,q \in[1,\infty]$ and $u \in L^p(\mathcal{G})\cap L^q(\mathcal{G})$} \,.
\end{equation}
The semigroup $\{T_p(t)\}_{t\ge0}$ is strongly continuous if $p \in[1,\infty)$ and holomorphic if $p \in(1,\infty)$.
\end{prop}


\subsection{Metric trees}\label{metre}

\begin{definition}\label{detree}
\noindent $(i)$ A metric graph $\mathcal{T}$ associated with a tree, with a {\em root} vertex $O$ singled out, is called a {\em metric tree}. The quantity $h(\mathcal{T}):=\sup_{x\in\mathcal{T}}\r(x)$\,, where $\r(x):=d(x,O)$ $(x\in \mathcal{T})$, is called the {\em height} of $\mathcal{T}$.

\noindent $(ii)$ The {\em generation} ${\rm gen}(v)$ of a vertex $v \in \mathcal{T}$ is the number of vertices which lie on the unique path connecting $v$ with $O$ (including the starting point but excluding the final point). The {\em branching number} $b(v)$ of $v \in \mathcal{T}$ is the number of edges emanating from $v$.

\noindent $(iii)$ A tree $\mathcal{T}$ is {\em regular}, if all vertices of the same generation have equal branching numbers, and all edges emanating from them are of the same length. Clearly,

\noindent $(a)$ in a regular tree all vertices of the same generation have the same distance from the root;

\noindent $(b)$ a regular tree is uniquely determined by two {\em generating sequences} $\{b_n\}\equiv\{b_n(\mathcal{T})\}\subseteq\{0\}\cup\N$,
$\{\r_n\}\equiv\{\r_n(\mathcal{T})\}\subseteq\R_+$, with $\{\r_n\}$ increasing and $\r_0=0$, such that
$$
n={\rm gen}(v)\,, \quad b_n=b(v)\,,\quad \r_n=\r(v) \quad\text{for any vertex $v$ of $\mathcal{T}$}
$$
(namely, $b_n$ is the common branching number of all vertices of generation $n$).

If the tree $\mathcal{T}$ is regular, its {\em branching function} $\b:[0,h(\mathcal{T}))\mapsto \N$ is defined as follows:
\begin{equation}\label{bratree}
\b(\r):={\rm card}\{x\in\mathcal{T}\,|\, \r(x)=\r \}\,=\, \sum_{n=1}^\infty ( b_0b_1\dots b_{n-1})\chi_{(\r_{n-1},\r_n]}(\r)\,.
\end{equation}
\end{definition}
Clearly, $\b$ is an increasing step function continuous from the left, with jump $\b(\r_n^+)-\b(\r_n^-)=b_0b_1\dots b_{n-1}(b_n-1)$ at any point $\r=\r_n$ $(n\in\N)$. We always suppose that $b_0 =1$, $b_n\ge2$ for any $n\in\N$, ${\rm gen}(v)<\infty$ for all $v \in \mathcal{T}$, and $h(\mathcal{T})=\infty$\,. In particular, there holds $\pa\mathcal{T}=O$ (i.e., the boundary of a metric tree consists of its root).

\begin{definition}\label{symm}
 A function $f:\mathcal{T}\mapsto \R$ is called {\em symmetric} if
$$
x, y \in \mathcal{T},\;\,   \r(x)=\r(y) \quad  \Rightarrow \quad f(x)=f(y)\,.
$$
\end{definition}
Clearly, $f$ is symmetric if and only if there exists $\tilde f:[0,h(\mathcal{T}))\mapsto \R$ such that
\begin{equation}\label{ftf}
f(x)=\tilde f(\r) \quad\text{for any $x\in \mathcal{T}$ with $\r(x)=\r$}\,.
\end{equation}

It is immediately proven (see \cite{NS2}) that:
\begin{equation}\label{leftf}
\sum_{x\in\mathcal{T}\!\!,\,\r(x)=\r} f(x)\,=\, \tilde f(\r)\,\b(\r) \quad\text{for a.e. $x\in\mathcal{T}\,.$}
\end{equation}\,



\end{document}